\documentclass[11pt]{article}

\usepackage{IEEEtrantools}
\usepackage{latexsym}
\usepackage{amsfonts}
\usepackage{amssymb}
\usepackage{amsmath,amsfonts,amsthm}
\usepackage{mathrsfs}
\usepackage{enumerate}

\usepackage{multirow}
\usepackage{tikz}
\usetikzlibrary{arrows,decorations.markings}
\usepackage{bookmark}
\usepackage{hyperref}
\hypersetup{
     colorlinks   = true,
     citecolor    = blue
}

\newcommand{\bbbt}{\mathbb{T}}
\newcommand{\scrt}{\mathscr{T}}

\newcommand{\be}{\begin{equation}}
\newcommand{\ee}{\end{equation}}
\newcommand{\bea}{\begin{eqnarray}}
\newcommand{\eea}{\end{eqnarray}}
\newcommand{\bean}{\begin{eqnarray*}}
\newcommand{\eean}{\end{eqnarray*}}
\newcommand{\brray}{\begin{array}}
\newcommand{\erray}{\end{array}}
\newcommand{\biearray}{\begin{IEEEarray}{rCl}}
\newcommand{\eiearray}{\end{IEEEarray}}

\newcommand{\newsection}[1]{\setcounter{equation}{0}
\setcounter{dfn}{0}
\section{#1}}

\newtheorem{dfn}{Definition}[section]
\newtheorem{thm}[dfn]{Theorem}
\newtheorem{lmma}[dfn]{Lemma}
\newtheorem{ppsn}[dfn]{Proposition}
\newtheorem{crlre}[dfn]{Corollary}
\newtheorem{xmpl}[dfn]{Example}
\newtheorem{rmrk}[dfn]{Remark}

\newcommand{\bdfn}{\begin{dfn}\rm}
\newcommand{\bthm}{\begin{thm}}
\newcommand{\blmma}{\begin{lmma}}
\newcommand{\bppsn}{\begin{ppsn}}
\newcommand{\bcrlre}{\begin{crlre}}
\newcommand{\bxmpl}{\begin{xmpl}}
\newcommand{\brmrk}{\begin{rmrk}\rm}

\newcommand{\edfn}{\end{dfn}}
\newcommand{\ethm}{\end{thm}}
\newcommand{\elmma}{\end{lmma}}
\newcommand{\eppsn}{\end{ppsn}}
\newcommand{\ecrlre}{\end{crlre}}
\newcommand{\exmpl}{\end{xmpl}}
\newcommand{\ermrk}{\end{rmrk}}

\newcommand{\bbc}{\mathbb{C}}
\newcommand{\bbz}{\mathbb{Z}}
\newcommand{\bbn}{\mathbb{N}}

\newcommand{\cla}{\mathcal{A}}

\newcommand{\clh}{\mathcal{H}}

\newcommand{\clk}{\mathcal{K}}

\def \bbt {\mbox{\boldmath $t$}}

\newcommand{\prf}{\noindent{\it Proof\/}: }

\newcommand{\id}{\mbox{id}}

\def \qed { \mbox{}\hfill
$\Box$\vspace{1ex}}

\begin{document}

\author{\sc{Satyajit Guin and Bipul Saurabh}}
\title{Representations and Classification of the compact quantum groups $U_q(2)$ for complex deformation parameters}
\maketitle


\begin{abstract} 
In this article, we obtain a complete list of inequivalent irreducible representations of the compact quantum group $U_q(2)$ for non-zero complex deformation parameters $q$, which are not roots of unity. The matrix coefficients of these representations are described in terms of the little $q$-Jacobi polynomials. The Haar state is shown to be faithful and an orthonormal basis of $L^2(U_q(2))$ is obtained. Thus, we have an explicit description of the Peter-Weyl decomposition of $U_q(2)$. As an application, we discuss the Fourier transform and establish the Plancherel formula. We also describe the decomposition of the tensor product of two irreducible representations into irreducible components. Finally, we classify the compact quantum groups $U_q(2)$.
\end{abstract}
\bigskip

{\bf AMS Subject Classification No.:} {\large 58}B{\large 32}, {\large 58}B{\large 34}, {\large46}L{\large 89}.

\bigskip


\newsection{Introduction}

In the theory of compact quantum group (CQG) developed by Woronowicz (\cite{Wor1},\cite{Wor2},\cite{Wor3}), the first nontrivial and the most accessible example is the $SU_q(2)$ for $q\in\mathbb{R}\setminus\{0\}$. It is widely studied in the literature through different perspective. Its representation theory (\cite{Koor-1989aa},\cite{Masuda-1991aa})  has striking similarity with that of its classical analogue $SU(2)$. Chakraborty-Pal \cite{CP-2003aa} studied its geometry along the line of Connes \cite{Con-1994aa}, and constructed a finitely summable $K$-homologically nontrivial spectral triple that is equivariant under its own comultiplication action. Connes proved regularity of this spectral triple and computed  its  local chern character \cite{Con-2004aa}. This builds a bridge  between the compact quantum groups of Woronowicz and Noncommutative geometry program of Connes. In a recent work \cite{KasMeyRoyWor-2016aa}, Woronowicz et al. defined a family of $q$-deformations of $SU(2)$ for $q\in\mathbb{C}\setminus\{0\}$. This agrees with the earlier defined $SU_q(2)$ when $q$ is real. Note that for $q\in\mathbb{C}\setminus\mathbb{R},\,SU_q(2)$ is not a CQG, but a braided compact quantum group in a suitable tensor category. In \cite{MeyRoyWor-2016aa}, it is shown that the quantum analogue of the semidirect product construction for groups turns the braided quantum group $SU_q(2)$ into a genuine CQG. This CQG is the coopposite of the compact quantum group $U_q(2)$ defined in \cite{ZhaZha-2005aa}. In this article, our object of study is this CQG $U_q(2),\,q\in\mathbb{C}\setminus\{0\}$, for generic parameter values (i,e. not a root of unity) with the intention to prepare the ground for investigating its noncommutative geometric aspects.

For real $q$, Koelink \cite{Koel-1991aa} initiated the study of $U_q(n)$. He begins with a deformation of the algebra of polynomials on $U(n)$. This algebra, denoted by $Pol(U_q(n))$, is a Hopf $\star$-algebra and is completed into a unital $C^*$-algebra $C(U_q(n))$ which yields a compact quantum group in the sense of Woronowicz. The matrix coefficients and Peter-Weyl decomposition is obtained  explicitly in \cite{MHM-1993aa}. The $n=2$ case is studied further in \cite{Wys-2004aa}, and it is shown that as a compact quantum group $U_q(2)$ is some twisted tensor product (closely related to Drinfeld's quantum double construction) of $SU_q(2)$ and $C(\mathbb{T})$. Note that at the level of $C^*$-algebras, $C(U_q(2))\cong C(SU_q(2))\otimes C(\mathbb{T})$ for real $q$. $U_q(2)$ has also been studied further from representation category viewpoint (e.g. \cite{Bic-1999aa},\cite{Mro-2014aa} and references therein). Now, we turn to the complex $q$ situation which is our object of study. The main difficulty in defining $U_q(2)$ for complex $q$, unlike the real case, is to establish a suitable $\star$-structure on the Hopf algebra $\mathcal{O}(GL_q(2))$. Zhang-Zhao \cite{ZhaZha-2005aa} suceeded in defining the CQG $U_q(2)$ in Woronowicz's framework following the FRT approach by taking Hayashi's $R$-matrix with $n=2$. This agrees with the earlier defined $U_q(2)$ when $q$ is real. The $\theta$-deformation $U_\theta(2)$ of the classical Lie group $U(2)$, constructed by Connes-Dubois-Violette in \cite{CV-2002aa}, is a special case of $U_q(2)$ for $q=exp(i\theta)$ \cite{Zha-2006aa}. For $|q|\neq 1,\,U_q(2)$ is a CQG {\it of non-Kac type} and for $|q|=1$ it is {\it of Kac type}. While for real $q$ there is a nondegenerate pairing between $U_q(\mathfrak{gl}(2))$ and the Hopf algebra $\mathcal{O}(U_q(2))$ associated with $U_q(2)$ (Page $316$ and $440$ in \cite{KliSch-1997aa}), for non-real complex numbers no such pairing has been established yet. Hence, it is not clear how $\mathcal{O}(U_q(2))$ can be viewed as certain subset of the dual of $U_q(\mathfrak{gl}(2))$. This makes the study of $U_q(2)$ more involved, but at the same time interesting and worth investigating. 
 
Our goal is to understand the geometry of $U_q(2)$ along the line of Connes, and this article is a first step towards that. For this purpose, the appropriate object to look for is an equivariant spectral triple (unbounded $\mathcal{K}$-cycle). It captures not only the geometric data of a CQG but also its comultiplication action. To construct spectral triples, one first needs to obtain explicit inequivalent irreducible corepresentations and then, the matrix coefficients need to be closely analysed. Here, the abstract realization of representation theory of $U_q(2)$ will not be of much use since, to show that the Dirac operator has bounded commutator with the matrix coefficients, an explicit description of them is needed. That is, even if one concludes that $U_q(2)$ has the same representation type of some classical Lie group, it does not help to produce an equivariant Dirac operator on $U_q(2)$. Faithfulness of the Haar state $h$ is desirable to get a faithful representation of the underlying $C^*$-algebra $C(U_q(2))$ on the GNS Hilbert space $L^2(h)$. Obtaining a precise orthonormal basis of $L^2(h)$ comes in handy in many situations, especially in determining action of the generators on these basis elements. Bounds on the ratio of the norm of the matrix coefficients chosen appropriately put a restriction on the growth of the singular values of an equivariant Dirac operator, which can be crucial in obtaining optimal spectral triple and computing the spectral dimension. Thus, it is also important to obtain norm of the matrix elements of irreducible corepresentations of $U_q(2)$. This article deals with these problems. 

Brief description of our work is the following. We separate three cases namely, $|q|<1, |q|=1$ and $|q|>1$, as the representation theory of the $C^*$-algebra $C(U_q(2))$ is different accordingly. First we deal the $|q|<1$ case. In Sec. \ref{Sec 2}, we establish a faithful representation of $C(U_q(2))$, and prove faithfulness of the Haar state. Note that in \cite{Wor2}, faithfulness of the Haar state is shown only on a dense subalgebra and therefore, this is not automatic on the $C^*$-algebra $C(U_q(2))$. In Sec. \ref{Sec 3}, using actions of $\mathbb{T}$ we put a $\mathbb{Z}^3$ grading on the Hopf $\star$-algebra $\mathcal{O}(U_q(2))$. In Sec. \ref{Sec 4}, we describe all finite dimensional inequivalent irreducible representations of $U_q(2)$, find their matrix coefficients explicitly, compute their norm, and obtain the Peter-Weyl decomposition. This gives an explicit orthonormal basis of $L^2(U_q(2))$ consisting of  normalized matrix elements. We have also expressed these matrix coefficients in terms of the little $q$-Jacobi polynomials. The Fourier transform on $U_q(2)$ is discussed, and the Plancherel formula is established. Then in Sec. \ref{Sec 5}, we give an explicit decomposition of the tensor product of two irreducible representations into irreducible components. Although it looks similar to that in the case of $SU_q(2)$, the proof is completely different. In the case of $SU_q(2)$, one first finds formula for the tensor product  decomposition of representations  of $U_q(\mathfrak{sl}_2)$ in the Drinfeld-Jimbo side (\cite{Drin-1988aa},\cite{Jim-1986aa}), and then through nondegenerate pairing  between $U_q(\mathfrak{sl}_2)$ and  the canonical Hopf algebra of $SU_q(2)$, one gets the result in the dual side.  Since, no such pairing is known yet for the case of $U_q(2)$ with $q\in\mathbb{C}\setminus\{0\}$, one can not follow this line of argument. Our approach uses explicit description of the matrix coefficients and $L^{\infty}$ functional calculus.

Next, in Sec. \ref{Sec 6}, we deal the case of $|q|=1$. We define faithful $C^*$-representation of $C(U_q(2))$ in this case, and prove faithfulness of the Haar state. Then, the Peter-Weyl decomposition and the Plancherel formula follows from similar computations done in the case of $|q|<1$, and hence we only mention the final results. This subsumes earlier investigations in (\cite{ZhaZha-2005aa},\cite{Zha-2006aa}). Finally, in Sec.  \ref{Sec 7}, we have classified the CQG $U_q(2)$ and it turns out that for nonzero complex numbers $q$ and $q^{'}$ which are not roots of unity, $U_q(2)$ and $U_{q^{'}}(2)$ are isomorphic if and only if $q^{'} \in\{q,\overline{q},\frac{1}{q},\frac{1}{\overline{q}}\}$. This also justifies that for $|q|\neq 1$, it is enough to do the computations for $|q|<1$ only.
\bigskip

\textbf{Notations~:} The following notations are used throughout the article.
\begin{enumerate}[{(i)}]
\item $\mathbb{A}_q$ denotes the $C^*$-algebra $C(U_q(2))$ and $\mathcal{A}_q$ denotes the underlying Hopf $\star$-algebra $\mathcal{O}(U_q(2))$.
\item $a,b,D$ will always denote the generators of $\mathbb{A}_q$.
\item $S$ will always denote the coinverse of $\mathcal{A}_q$.
\item $t^l_{ij}D^k$ denotes the matrix coefficients of irreducible corepresentations of $\mathcal{A}_q$.
\item `$i$' will always denote an element in $\frac{1}{2}\bbn$, and whenever the complex inderminate appears we write it $\sqrt{-1}$.
\item The standard orthornormal basis of $\ell^2(\bbn)$ (or $\ell^2(\bbz)$) is denoted by $\{e_n:n\in\bbn\}$ (or $n\in\bbz$ depending on context).
\item $N$ denotes the number operator $e_n\mapsto ne_n$ acting on $\ell^2(\bbn)$ or $\ell^2(\bbz)$ depending on context.
\item $V$ denotes the right shift operator $e_n\mapsto e_{n+1}$ acting on $\ell^2(\bbn)$, and $U$ denotes the unitary shift operator $e_n\mapsto e_{n+1}$ acting on $\ell^2(\bbz)$.
\end{enumerate}
\bigskip


\newsection{Faithful representation and the Haar state}\label{Sec 2}

Throughout the article, $q$ denotes a non-zero complex number which is not a root of unity, and we use the notations $\mathbb{C}^*:=\mathbb{C}\setminus\{0\},\,\mathbb{R}^*:=\mathbb{R}\setminus\{0\}$. We first recall the compact quantum group $U_q(2)$ as defined in \cite{ZhaZha-2005aa}. The $C^*$-algebra $C(U_q(2))$, to be denoted by $\mathbb{A}_q$, is the universal $C^*$-algebra generated by $a,b,D$ satisfying the following relations~:
\begin{IEEEeqnarray}{lCl} \label{relations}
ba&=& q ab, \qquad a^*b=qba^*, \qquad   \qquad \qquad bb^*=b^*b, \qquad \qquad aa^*+bb^*=1, \nonumber \\
aD&=&Da, \qquad bD=q^2|q|^{-2}Db,  \qquad DD^*=D^*D=1,  \qquad a^*a+|q|^2b^*b=1.  
\end{IEEEeqnarray}
The compact quantum group structure is given by the comultiplication $\Delta:\mathbb{A}_q\longrightarrow\mathbb{A}_q\otimes\mathbb{A}_q$ defined as follows~:
\begin{IEEEeqnarray}{lCl} \label{comul}
 \Delta(a)=a \otimes a-\bar{q}b \otimes Db^*\quad,\quad\Delta(b)=a\otimes b+b\otimes Da^*\quad,\quad\Delta(D)=D \otimes D. 
\end{IEEEeqnarray}
Let $\mathcal{O}(U_q(2))$ be the $\star$-subalgebra of $C(U_q(2))$ generated by $a,b$ and $D$. We will denote it by $\cla_q$. 
The Hopf $\star$-algebra structure on it is given by the following~: 
\begin{IEEEeqnarray*}{lCl}
\mbox{antipode:} \quad S(a)=a^*,\,\,S(b)=-qbD^*,\,\,S(D)=D^*,\,\,S(a^*)=a,\,\,S(b^*)=-(\bar{q})^{-1}b^*D\,,\\ 
\mbox{counit:} \qquad\epsilon(a)=1,\,\,\epsilon(b)=0,\,\,\epsilon(D)=1\,.
\end{IEEEeqnarray*}
We divide the article in two cases namely, the case of $|q|=1$ and $|q|\neq 1$. First we restrict our attention to the case of $|q|<1$ and $q\neq 0$ till Section $5$. In Section $7$, where we prove the classification of $U_q(2)$, we will see that for $|q|\neq 1$, it is enough to provide computations only for the case of $|q|<1$. The case of $|q|=1$ is dealt separately in Section $6$.
\medskip

Fix any $q\in\mathbb{C}^*$ with $|q|<1$ and let $\theta=\frac{1}{\pi}\arg{(q)}$. The $C^*$-algebra $\mathbb{A}_q=C(U_q(2))$ can be realized more concretely as follows. Let $\clh$ be the Hilbert space $\ell^2(\bbn)\otimes\ell^2(\bbz)\otimes\ell^2(\bbz)$. Consider the right shift operator $V:e_n\mapsto e_{n+1}$ acting on
$\ell^2(\bbn)$ and the bilateral shift $U:e_n\mapsto e_{n+1}$ acting on $\ell^2(\bbz)$. Define the following representation $\pi$ of $\mathbb{A}_q$ on $\clh$~:
\begin{IEEEeqnarray}{rCl}\label{representation}
\pi(a)= \sqrt{1-|q|^{2N}}\,V\otimes I\otimes I, \quad 
\pi(b) =q^N \otimes U\otimes I, \quad
\pi(D) = I\otimes e^{-2 \pi\sqrt{-1}\theta N}\otimes U.
\end{IEEEeqnarray}
Let $\scrt$ be the Toeplitz algebra generated by the right shift operator $V$ acting on $\ell^2(\bbn)$. Denote by $\,\bbt:t\longmapsto t$ the identity map on $\bbbt$. Then, 
the homomorphism $\,\bbt\longmapsto U$ is an isomorphism between $C(\bbbt)$ and the $C^*$-subalgebra of 
$\mathcal{B}(\ell^2(\bbz))$ generated by the bilateral shift $U$ acting on $\ell^2(\bbz)$. Let $\sigma:\scrt\longrightarrow C(\bbbt)$ be the symbol map sending $V$ to $\bbt$, and $ev_t: C(\bbbt)\longrightarrow \bbc$ denotes the evaluation map at $t$.
\bppsn
The representation $\pi$ of $\mathbb{A}_q$ defined above is faithful.
\eppsn
\prf Since $\sqrt{1-|q|^{2N}}\,V,\,q^N\in\scrt$, one can view $\mathbb{A}_q\subseteq \scrt \otimes\mathcal{B}(\ell^2(\bbz)) \otimes C(\bbbt)$. 
Consider the following homomorphisms,
\begin{IEEEeqnarray*}{rCl}
\sigma\otimes I\otimes I &:& \scrt\otimes\mathcal{B}(\ell^2(\bbz))\otimes C(\bbbt)\longrightarrow C(\bbbt)\otimes\mathcal{B}(\ell^2(\bbz)) \otimes C(\bbbt)\,,\\
I\otimes I\otimes ev_t &:& \scrt\otimes\mathcal{B}(\ell^2(\bbz))\otimes C(\bbbt)\longrightarrow C(\bbbt)\otimes\mathcal{B}(\ell^2(\bbz))\,.
\end{IEEEeqnarray*} 
Let  $\varphi$ and $\psi$ be the restriction of $\sigma \otimes I \otimes I$ and $I \otimes I \otimes ev_t$ respectively to $\mathbb{A}_q\,$. One has 
\[
 \varphi(\pi(a))=\bbt \otimes I \otimes I\quad,\quad\varphi(\pi(b))= 0 \quad,\quad\varphi(\pi(D))=I \otimes e^{-2 \pi \sqrt{-1}\theta N }\otimes \bbt\,\,;
\]
and 
\[
 \psi(\pi(a))=\sqrt{1-|q|^{2N}}\,V\otimes I\quad,\quad\psi(\pi(b))= q^N\otimes U\quad,\quad\psi(\pi(D))=tI \otimes e^{-2 \pi\sqrt{-1}\theta N }\,\,.
\]
It is not difficult to see that all one dimensional irreducible representations of $\mathbb{A}_q$ factor through the homomorphism $\varphi\circ\pi$, and all infinite dimensional representations of $\mathbb{A}_q$ described in Lemma $3.1$ and Thm. $3.2$ in \cite{ZhaZha-2005aa} factor through the homomorphism $\psi\circ\pi$. This proves that any irreducible representation of $\mathbb{A}_q$ factors through the map $\pi$. Hence, image of an element in $\ker(\pi)$ under any irreducible representation of $\mathbb{A}_q$ is zero, which implies that the element is zero because norm of any element in
a $C^*$-algebra is the supremum of the norms of the image of that element under all irreducible representations. This proves the claim.\qed
\brmrk
Note that Lemma $3.1$ in \cite{ZhaZha-2005aa} is for the case of $|q|>1$. In the case of $|q|<1$, one needs to define $H_0$ as $\ker(a^*)$.
\ermrk
From now on, thanks to the above proposition, we will identify $a$, $b$ and $D$ with $\pi(a)$, $\pi(b)$ and $\pi(D)$ 
respectively. Moreover, we view $\mathbb{A}_q$ as the $C^*$-subalgebra of $\mathcal{B}\big(\ell^2(\bbn)\otimes\ell^2(\bbz)\otimes\ell^2(\bbz)\big)$ generated by $\pi(a)$, $\pi(b)$ and $\pi(D)$. For $n,l \in \bbz$ and $ m,k \in \bbn$, define
\begin{IEEEeqnarray*}{rCl}
 \langle n,m,k,l \rangle = \begin{cases}
                            a^n b^m(b^*)^k D^l & \mbox{ if } n \geq 0, \cr
                            (a^*)^{-n} b^m(b^*)^k D^l & \mbox{ if } n \leq 0. \cr
                           \end{cases}
\end{IEEEeqnarray*} 
\bthm[\cite{ZhaZha-2005aa}]\label{basis}
The set $\{ \langle n,m,k,l \rangle: n,l \in \bbz, m,k \in \bbn\}$ forms a linear basis of $\mathcal{O}(U_q(2))$ for all $q\in\mathbb{C}^*$.
\ethm
\bthm[\cite{ZhaZha-2005aa}]\label{haar} 
The Haar state $h:C(U_q(2))\longrightarrow\bbc$ is given by the following,
\begin{IEEEeqnarray*}{rCl}
h(x)= (1-|q|^2)\sum_{n=0}^{\infty} |q|^{2n} \langle e_{n,0,0}\,,\,\pi(x)e_{n,0,0}\rangle\,,
\end{IEEEeqnarray*}
where $\{e_{n,r,s}\}$ denotes the standard orthonormal basis of $\ell^2(\bbn)\otimes\ell^2(\bbz)\otimes\ell^2(\bbz)$. Moreover,  one has
\begin{IEEEeqnarray*}{rCl}
 h( \langle n,m,k,l \rangle)= \begin{cases}
                               \frac{1-|q|^2}{1-|q|^{2(m+1)}} & \mbox{ if } m=k, \mbox{ and } n=l=0, \cr
                               0 & \mbox{ otherwise .}
                              \end{cases}
\end{IEEEeqnarray*}
\ethm
In this case of $|q|<1$, the Haar state is not a trace as
\begin{center}
$h(a^*a-aa^*)=(1-|q|^2)h(bb^*)=\frac{1-|q|^2}{1+|q|^2}\neq 0\,.$
\end{center}
Let $\Delta_{\bbbt}$ denotes the coproduct on $C(\bbbt)$ and $\phi:\mathbb{A}_q\longrightarrow C(\bbbt)$ be the homomorphism given by $\,\phi(a)=1,\,\phi(b)=0$ and $\phi(D)=\bbt$. One can check that $\Delta_{\bbbt} \circ \phi=(\phi \otimes \phi)\circ \Delta\,$. Thus, $\bbbt$ is a quantum subgroup of $U_q(2)$. Hence, $\bbbt$ acts on $\mathbb{A}_q$ by the formula $\Phi(x)=(\mbox{id}\otimes\phi)\Delta\,$. In such a case, one defines the quotient space $U_q(2)/\bbbt$ as follows,
\[
 C(U_q(2)/\bbbt)=\{x\in\mathbb{A}_q:(\mbox{id}\otimes\phi)\Delta(x)=x\otimes 1\}.
\]
The conditional expectation $E:\mathbb{A}_q\longrightarrow C(U_q(2)/\bbbt)$ is defined as $(\mbox{id}\otimes(h_{\bbbt}\circ\phi))\circ\Delta\,,$ where $h_{\bbbt}$ denotes the Haar state on $\bbbt$. 
\blmma The $C^*$-algebra $C(U_q(2)/\bbbt)$ is the $C^*$-subalgebra of $\mathbb{A}_q$ generated by $a$ and $b$.
\elmma
\prf Using the continuity of $E$ and by the open mapping theorem, one has
\begin{IEEEeqnarray}{rCl} \label{eq1}
C(U_q(2)/\bbbt)=E(\mathbb{A}_q)=E(\overline{\cla_q})=\overline{E(\cla_q)}\,.
\end{IEEEeqnarray}
Applying the formula $E=(\mbox{id}\otimes(h_{\bbbt}\circ\phi))\circ\Delta\,$, one can see that
\[
E(a)=a\,,\,E(a^*)=a^*\,,\,E(b)=b\,,\,E(b^*)=b^*\,\,\mbox{and}\,\,E(D^m)=0\,\,\mbox{for}\,\,m\in\bbz\setminus\{0\}.
\]
This proves that $C(U_q(2)/\bbbt)$ contains the $\star$-subalgebra generated by $a$ and $b$. Using the property $E(xy)=xE(y)$ for $x\in C(U_q(2)/\bbbt)$ and $y\in\mathbb{A}_q\,$, we get the following,
\[
E(\langle n,m,k,l \rangle)=\langle n,m,k,0\rangle E(D^l)=\begin{cases}
                                                          \langle n,m,k,0\rangle & \mbox{ if } l \neq 0,\cr
                                                          0 & \mbox{ if } l = 0.\cr
                                                         \end{cases}
\]
Hence, $E(\cla_q)$ is the $\star$-algebra generated by $a$ and $b$. This, along with the eqn. \ref{eq1}, proves the claim.
\qed

Denote by $\mathbb{A}_q^0$ the $C^*$-subalgebra $C(U_q(2)/\bbbt)$. Let $ C(SU_{|q|}(2))$ be the $C^*$-subalgebra of $\mathcal{B}(\ell^2(\bbn)\otimes\ell^2(\bbz)\otimes\ell^2(\bbz))$ generated by the operators $\sqrt{1-|q|^{2N}}\,V\otimes I\otimes I$ and $|q|^N\otimes U\otimes I$. 
\blmma\label{equal}
For $|q|<1$, one has $\mathbb{A}_q^0=C(SU_{|q|}(2))$.
\elmma
\prf To see $\,C(SU_{|q|}(2))\subseteq\mathbb{A}_q^0\,$, first observe that $b^*b=|q|^{2N}\otimes I\otimes I$. Using the spectral decomposition of $\,b^*b$, one can see that $\,p_i\otimes I\otimes I\in\mathbb{A}_q^0\,\,\forall\,i\in\bbn$, where $p_i$ is the rank one projection acting on $\ell^2(\bbn)$ by $p_i(e_n)=\delta_{in}e_n$. Hence, $T_n=\big(\sum_{i=0}^n\,\frac{|q|^i}{q^i}p_i\otimes I\otimes I\big)b\in\mathbb{A}_q^0\,$. Then, we have
\[
 \|T_n-|q|^N\otimes U\otimes I\|\leq |q|^{n+1}\|U\|\|I\|=|q|^{n+1}.
\]
This shows that $\lim_{n\rightarrow\infty}T_n =|q|^N\otimes U\otimes I$ and hence, $|q|^N\otimes U\otimes I \in \mathbb{A}_q^0\,$. The reverse inclusion follows from similar argument. 
\qed
\brmrk
The above Lemma and Lemma $3.2$ of \cite{PalSun-2010aa} proves that $\mathbb{A}_{q_1}^0=\mathbb{A}_{q_2}^0$ for $q_1,q_2\in\mathbb{C}^*$ with $|q_1|,|q_2|<1$.  It is not difficult to see that $\mathbb{A}_q^0$ is 
same as $C(SU_q(2))$ mentioned in \cite{KasMeyRoyWor-2016aa} and therefore, the isomorphism of $\mathbb{A}_q^0$ follows from Thm. $2.3$ in \cite{KasMeyRoyWor-2016aa}. But here we are saying slightly more. Lemma \ref{equal} says that the identify map is an isomorphism between $\mathbb{A}_q^0$ and $\mathbb{A}_{|q|}^0$. This is crucial in getting the faithfulness of the Haar state.
\ermrk
\bthm\label{faithful}
The Haar state $h$ on the quantum group $U_q(2)$ is faithful.
\ethm
\prf To prove the claim, we will apply Lemma $2.1$ and Propn. $2.2$ of \cite{Nag-1993aa}. Faithfulness of $h_{\bbbt}$ is well-known. Moreover, the counit $\epsilon_{\bbbt}\,$, initially defined on the algebra of polynomials in $\bbt$ and $\bbt^{-1}$,  is same as evaluation at $1$ and hence, it can be extended to $C(\bbbt)$. It follows from Thm. \ref{haar} and Lemma \ref{equal} that the restriction $\left.h\right|_{\mathbb{A}_q^0}$ of the  Haar state $h$ of $\mathbb{A}_q$ to $\mathbb{A}_q^0$ is same as the Haar state of $C(SU_{|q|}(2))$. Hence, it follows from Thm. $1.1$ in \cite{Nag-1993aa} that $\left.h\right|_{\mathbb{A}_q^0}$ is a faithful state on $\mathbb{A}_q^0\,$. Combining all these facts and using Lemma $2.1$ of \cite{Nag-1993aa}, we get the claim.\qed
\bigskip


\newsection{Graded decomposition of \texorpdfstring{$\mathcal{O}(U_q(2))$}{}}\label{Sec 3}

In this section we decompose $\cla_q=\mathcal{O}(U_q(2))$ in a way similar to what is done in (page 105, \cite{KliSch-1997aa}) for $\mathcal{O}(SL_q(2))$, but  with respect to three different coactions. We assign three integers to each basis element and  put a $\bbz^3$ grading on $\cla_q$. Consider the Hopf $\star$-algebra $\mathcal{O}(\bbbt)=\bbc[z,z^{-1}]$ where,
\begin{center}
$z^*=z^{-1}\quad,\quad\Delta(z)=z\otimes z\quad,\quad S(z)=z^{-1}\quad,\quad\epsilon(z)=1\,.$
\end{center}
Define a $\star$-homomorphism $\phi_\bbbt:\cla_q \longrightarrow \mathcal{O}(\bbbt)$ given by $\phi_\bbbt(a)=z$, $\phi_\bbbt(b)=0$ and $\phi_\bbbt(D)=1$. Using this homomorphism, one can put a natural left $\mathcal{O}(\bbbt)$ comodule structure $L_\bbbt:  \cla_q\longrightarrow\mathcal{O}(\bbbt) \otimes \cla_q$ and a right $\mathcal{O}(\bbbt)$ comodule structure  $R_\bbbt: \cla_q \longrightarrow \cla_q\otimes \mathcal{O}(\bbbt)$ defined as follows;
\begin{center}
$L_\bbbt=(\phi_\bbbt\otimes\id)\circ\Delta\qquad\mbox{and}\qquad R_\bbbt=(\id\otimes\phi_\bbbt)\circ\Delta\,.$
\end{center}
For $m,n \in \bbz$, define
\begin{center}
$\cla_q[m,n]=\{x \in \cla_q: L_\bbbt(x)= z^m \otimes x \mbox{ and } R_\bbbt(x)= x \otimes z^n\}\,.$
\end{center}
Let $\zeta=bb^*$.  If $m-n$ is even, define
\begin{IEEEeqnarray}{rCl}\label{defn of e_m,n}
e_{m,n}=\begin{cases}
          a^{\frac{m+n}{2}} b^{\frac{m-n}{2}} & \mbox{ if }  m+n\geq 0, m\geq n, \cr
           a^{\frac{m+n}{2}} (b^*)^{\frac{n-m}{2}} & \mbox{ if }  m+n\geq 0, m\leq n, \cr
           b^{\frac{m-n}{2}} (a^*)^{\frac{-m-n}{2}} & \mbox{ if }  m+n\leq 0, m\geq n, \cr
          (b^*)^{\frac{n-m}{2}} (a^*)^{\frac{-m-n}{2}} & \mbox{ if }  m+n\leq 0, m\leq n. \cr
          \end{cases} 
 \end{IEEEeqnarray}        

\bppsn\label{cla}
One has the followings.
\begin{enumerate}[(i)]
\item $\cla_q[m,n]\cla_q[p,q] \subset \cla_q[m+p,n+q]$.
\item $\langle n,m,k,l \rangle \in \cla_q[n+m-k, n-m+k]$.
\item $\cla_q[0,0]$ is a commutative $\star$-subalgebra of $\cla_q$ generated by $D$, $D^*$ and $bb^*$.
\item $\cla_q[m,n]$ is a $\cla_q[0,0]$-bimodule. 
\item $\cla_q= \oplus_{m,n \in \bbz} \cla_q[m,n]$.
\item If $m-n$ is odd, then $\cla_q[m,n]=\{0\}$.
\item If $m-n$ is even, then $\cla_q[m,n]=e_{m,n}\cla_q[0,0]=\cla_q[0,0]e_{m,n}\,$. 
\end{enumerate}
\eppsn
\prf It follows from the definition of $\cla_q[m,n]$ and the fact that $\{\langle n,m,k,l\rangle:n,l\in\bbz,m,k\in\bbn\}$ forms a linear basis of $\cla_q$.\qed

From part $(i)$ of the above Lemma, we see that this decomposition of $\cla_q$ into $\cla_q[m,n]$'s is a $\bbz^2$ grading on $\cla_q$. To put a $\bbz^3$ grading, consider another $\star$-homomorphism $\psi_\bbbt:\cla_q \longrightarrow \mathcal{O}(\bbbt)$ given by $\psi_\bbbt(a)=1,\,\psi_\bbbt(b)=0$ and $\psi_\bbbt(D)=z$. Define the coaction $L_\bbbt^{\psi}=(\psi_\bbbt\otimes\id)\circ\Delta\,$. Let
\begin{center}
$\cla_q[m,n,r]=\{x \in \cla_q[m,n]: L_\bbbt^{\psi}(x)= z^r \otimes x \}\,.$
\end{center}
\bppsn\label{clafurther}
One has the followings.
\begin{enumerate}[(i)]
\item $\cla_q[m,n,r]\cla_q[p,q,s] \subset \cla_q[m+p,n+q, r+s]$.
\item $\cla_q[m,n,r]^*=\cla_q[-m,-n,-r]$.
 \item $\langle n,m,k,l \rangle \in \cla_q[n+m-k, n-m+k,l]$.
\item $\cla_q[0,0,0]$ is isomorphic to 
$\bbc[\zeta]$.
 \item $\cla_q[m,n,r]=D^r\cla_q[m,n,0]=\cla_q[m,n,0]D^r$ for $r \in \bbz$.
 \item If $m-n$ is even, $\cla_q[m,n]=\oplus_{r \in \bbz}\cla_q[m,n,r]=\oplus_{r \in \bbz} D^r \cla_q[m,n,0]=\oplus_{r \in \bbz} D^r e_{m,n}\bbc[\zeta]$.
 \item $S(\cla_q[m,n,r])=\cla_q[-n,-m,-r]$.
\end{enumerate}
\eppsn
\prf Part $(i)$ and $(ii)$ follow from the fact that  $L_\bbbt$, $R_\bbbt$ and  $L_\bbbt^{\psi}$ are $\star$-preserving homomorphisms. Part $(iii)$ follows immediately from the definition. Observe that $D^r \in \cla_q[0,0,r]$ for $r \in \bbz$. Using this and part $(iii)$ of Propn. $\ref{cla}$, we get part $(iv)$. Part $(v)$ follows directly from part $(i)$. Part $(vi)$ follows from part $(iii)$, part $(vii)$ of Propn. $\ref{cla}$, and the fact that $\{ \langle n,m,k,l \rangle: n,l \in \bbz, m,k \in \bbn\}$ forms a linear basis of $\cla_q$. For part $(vii)$, observe that 
\begin{center}
$S(\zeta)=\zeta\quad,\quad S(D)=D^*\quad,\quad S(D^*)=D\quad\mbox{and}\quad S(e_{m,n})=Ce_{-n,-m}\,.$
\end{center}
where, $C$ is a nonzero constant. Combining this with part $(v)$, we get that
\begin{center}
$S(\cla_q[m,n,r])\subseteq\cla_q[-n,-m,-r]\,.$
\end{center}
To get the reverse inclusion, take any $y\in\cla_q[-n,-m,-r]$. From part $(ii)$, we have 
\begin{center}
$y^*\in \cla_q[n,m,r]\Rightarrow S(y^*)\in\cla_q[-m,-n,-r]\Rightarrow S(y^*)^*\in\cla_q[m,n,r]\,.$ 
\end{center}
Since $S(S(y^*)^*)=y$, we get the claim. 
\qed
\bppsn \label{deg same} If $m,n,r \in \bbz$ with $m-n$ is even, $\cla_q[m,n,r]$ is a free $\bbc[\zeta]$-left (or right)  module with basis $e_{m,n}D^r$. 
Moreover, if $f(\zeta)e_{m,n}D^r=e_{m,n}D^r g(\zeta)=e_{m,n}h(\zeta)D^r$  for polynomials $f,g,h$ then degree of $f,g$ and $h$ are same.  
\eppsn 
\prf By part $(v)$ of Propn. $\ref{clafurther}$, one can see that $e_{m,n}D^r$ generates  $\cla_q[m,n,r]$. 
Using the commutation relations in \ref{relations} and Thm. \ref{basis}, it follows that for any polynomial $P$, $P(\zeta)e_{m,n}D^r=0$ implies $P(\zeta)=0$. This proves that  
$e_{m,n}D^r$ is a $\bbc[\zeta]$-basis of $\cla_q[m,n,r]$. 
Since  $\zeta^ke_{m,n}=Ce_{m,n}\zeta^k$ for 
some nonzero constant, we get the last part of the claim.\qed
\bppsn\label{clafurtherfurther}
If $x\in\cla_q[m,n,r]$ and $(m,n,r)\neq (0,0,0)$, then $h(x)=0$.
\eppsn
\prf
For $z \in \bbbt$, consider the linear functionals $\beta_z=ev_z\circ \phi_{\bbbt}$ and $\gamma_z=ev_z\circ \psi_{\bbbt}$ on $\mathcal{A}_q$, where $ev_z$ is the evaluation map at $z$. Take any $x \in  \cla_q[m,n,r]$.  Using the invariance property of Haar measure, we have 
\[
h(x)=(\beta_z\otimes h)\circ \Delta (x)= (ev_z \otimes h) (\phi_{\bbbt}\otimes\mbox{id})\circ \Delta (x)= 
 (ev_z \otimes h)L_{\bbbt}(x)=(ev_z \otimes h)(z^m\otimes x)=z^mh(x)\,.
\]
for all $z \in \bbbt$. Similarly, considering  $h(x)=(h\otimes\beta_z)\circ\Delta (x)$ and $h(x)=(\gamma_z\otimes h)\circ \Delta (x)$, we get that
\[
 h(x)=z^nh(x) \qquad and \qquad h(x)=z^rh(x)\,.
\]
for all $z \in \bbbt$. Hence if $(m,n,r)\neq (0,0,0)$, we have $h(x)=0$.\qed
\brmrk
It follows that for any $x\in\cla_q,\,h(x^*)=\overline{h(x)}$ and $h(S(x))=h(x)$. To see this, take any $x\in\cla_q[m,n,r]$. From part $(ii)$ of Propn. $\ref{clafurther}$, it follows that $x^*\in\cla_q[-m,-n.-r]$. Hence, if $(m,n,r)\neq (0,0,0)$, then we have $h(x^*)=0=\overline{h(x)}$. If $x\in\cla_q[0,0,0]$, then by part $(iv)$ of Propn. $\ref{clafurther}$ we have 
$x=p(\zeta)$ where, $p$ is a polynomial with complex coefficients. Since $\zeta^*=\zeta$, we get that $h(x^*)=\overline{h(x)}$. The other part follows from a similar argument.
\ermrk
\bigskip


\newsection{The Peter-Weyl decomposition}\label{Sec 4}

In this section, we describe all the irreducible representations of $U_q(2)$ and obtain an orthonormal basis of $L^2(U_q(2),h)$ in terms of the matrix coefficients of these representations. Let us fix some notations. 
For $l\in\frac{1}{2}\bbn,\,m,n\in\bbn,\,\alpha\in\bbc$ and $q\in\bbc^*$, let
\begin{IEEEeqnarray*}{rCl}
c=-\bar{q}Db^* &\quad,\quad& d=Da^*\,.\\
I_l &=& \{-l, -l+1, \cdots , l-1,l\}\,.\\ 
(\alpha,q)_n &=& \begin{cases}
           1 & \mbox{ if } n=0\,,\cr
           \prod_{r=0}^{n-1}(1-\alpha q^r) & \mbox{ if } n>0\,.\cr
          \end{cases}\\
{n \choose m}_{q} &=& \frac{(q,q)_n}{(q,q)_m(q,q)_{n-m}}\,.\\
|m+1|_{|q|} &=& \sum_{k=0}^{m}|q|^{m-2k}\,. 
\end{IEEEeqnarray*}
Thus, $|q|^{m}|m+1|_{|q|}=\sum_{k=0}^{m}|q|^{2(m-k)}=\frac{1-|q|^{2(m+1)}}{1-|q|^2}\,$.

\subsection{Coaction of \texorpdfstring{$\mathcal{A}_q$}{} on the quantum plane}

Let $\mathcal{O}(\bbc_q^2)$ be the complex associative  unital  algebra generated by two symbols  
$x$ and $y$ satisfying $xy=qyx$. This is called the quantum plane.
\bppsn\label{rightaction}
The association 
\begin{IEEEeqnarray}{rCl} \label{actionformula}
x\longmapsto x\otimes Da^*+y\otimes b\quad\mbox{and}\quad y\longmapsto x\otimes(-\bar{q}Db^*)+y\otimes a
\end{IEEEeqnarray}
extends to a unique homomorphism $\psi_R: \mathcal{O}(\bbc_q^2)\longrightarrow\mathcal{O}(\bbc_q^2)\otimes\cla_q\,$. Moreover, $\psi_R$ is a right coaction of $\cla_q$ on $\mathcal{O}(\bbc_q^2)$. 
\eppsn
\prf It is enough to show that image of $xy$ and $qyx$ under the map $\psi_R$ are same. Using the commutation relations, we get that 
\begin{IEEEeqnarray*}{rCl}
 \psi_R(xy)&=&(x \otimes Da^*+y\otimes b)(x\otimes(-\bar{q}Db^*)+y \otimes a)\\
 &=& -\bar{q}x^2\otimes Da^*Db^*+y^2 \otimes ba+xy \otimes Da^*a-\bar{q}yx \otimes bDb^*\\
 &=&-\bar{q}x^2\otimes D^2a^*b^*+qy^2 \otimes ab+xy \otimes D(aa^*+(1-|q|^2)b^*b)-\frac{\bar{q}q^2}{q|q|^2}xy \otimes Dbb^*\\
 &=&-\bar{q}^2x^2\otimes D^2b^*a^*+qy^2 \otimes ab+xy \otimes Daa^*-|q|^2 xy \otimes Db^*b\\
 &=&-q\bar{q}x^2\otimes Db^*Da^*+qy^2 \otimes ab+qyx \otimes Daa^*-q\bar{q}xy \otimes Db^*b\\
 &=& q(x \otimes(-\bar{q}Db^*)+y\otimes a)(x \otimes Da^*+y \otimes b)\\
 &=& \psi_R(qyx)\,.
\end{IEEEeqnarray*}
To show that $\psi_R$ is a right coaction of $\cla_q$, one has to check the conditions $(\psi_R\otimes\mbox{id})\circ \psi_R=(\mbox{id}\otimes\Delta)\circ\psi_R$ and $(\mbox{id}\otimes\epsilon)\circ\psi_R=\mbox{id}$ on the generators $x$ and $y$ of $\mathcal{O}(\bbc_q^2)$, and this follows easily from the definition of $\psi_R$.
\qed
\bppsn
The association 
\begin{IEEEeqnarray}{rCl} \label{leftactionformula}
x\longmapsto Da^*\otimes x+(-\bar{q}Db^*)\otimes y\quad\mbox{and}\quad y\longmapsto b\otimes x+a\otimes y 
\end{IEEEeqnarray}
extends to a unique homomorphism $\psi_L:\mathcal{O}(\bbc_{\bar{q}}^2)\longrightarrow\cla_q\otimes\mathcal{O}(\bbc_{\bar{q}}^2)$. Moreover, $\psi_L$ is a left coaction of $\cla_q$ on $\mathcal{O}(\bbc_{\bar{q}}^2)$. 
\eppsn
\prf Proof is similar to that of Propn. \ref{rightaction}.\qed 

Let  $\mathcal{O}(\bbc_q^2)_{2l}$ be the vector subspace of $\mathcal{O}(\bbc_q^2)$ spanned by degree $2l$ homogeneous monomials in $x$ and $y$, i.e. $ \mathcal{O}(\bbc_q^2)_{2l}=\bigoplus_{j=-l}^{l}\bbc y^{l-j}x^{l+j}$.
It follows from \ref{actionformula} that $\mathcal{O}(\bbc_q^2)_{2l}$ is invariant under the map $\psi_R$. Hence, the restriction of the coaction $\psi_R$ to $ \mathcal{O}(\bbc_q^2)_{2l}$ gives a $(2l+1)$-dimensional
corepresentation of $\cla_q$ which will be denoted by $T_{l}$. In what follows, we will embed $ \mathcal{O}(\bbc_q^2)_{2l}$ in $\cla_q$ and find a corepresentation of $\cla_q$ on a vector subspace of $\cla_q$ which is equivalent to $T_l$. Let $\cla_q(a,b)$  be the subalgebra of $\cla_q$ generated by $a$ and $b$. The algebra $\cla_q(a,b)$ has a canonical $\cla_q$-comodule structure given by restriction of the  comultiplication $\Delta$ to $\cla_q(a,b)$. Since $ba=qab$, we get a homomorphism $\vartheta_R :\mathcal{O}(\bbc_q^2)\longrightarrow \cla_q(a,b)$ such that $\vartheta_R(x)=b$ and $\vartheta_R(y)=a$.
\bppsn \label{equi}
The homomorphism $\vartheta_R:\mathcal{O}(\bbc_q^2)\longrightarrow \cla_q(a,b)$ is an isomorphism between the right $\cla_q$-comodule algebra  $\mathcal{O}(\bbc_q^2)$ with the right coaction $\psi_R$ and the right $\cla_q$-comodule algebra  $\cla_q(a,b)$ with the comultiplication $\Delta$.
\eppsn
\prf Since $\{x^my^n:m,n\in\bbn\}$ and $\{b^ma^n:m,n\in\bbn\}$ are linear basis of $\mathcal{O}(\bbc_q^2)$ and $\cla_q(a,b)$ respectively, the homomorphism  $\vartheta_R$ is an algebra isomorphism. To prove that $\vartheta_R$ is an $\cla_q$-comodule isomorphism, it is enough to check the condition $(\vartheta_R \otimes \mbox{ id }) \circ \psi_R=\Delta \circ \vartheta_R$ on the generators $x$ and $y$ of $\mathcal{O}(\bbc_q^2)$. Applying the formula of comultiplication of $a,b$ and the coaction $\psi_R$, this follows easily.\qed

Similarly, one can define a homomorphism $\vartheta_L :\mathcal{O}(\bbc_{\bar{q}}^2)\longrightarrow\cla_q(a,c)$ such that $\vartheta_L(x)=c$ and $\vartheta_L(y)=a$, since $ca=\bar{q}ac$. This gives an isomorphism between the left $\cla_q$-comodule algebra  $\mathcal{O}(\bbc_{\bar{q}}^2)$ with the left coaction $\psi_L$ and the left $\cla_q$-comodule algebra  $\cla_q(a,c)$ with the comultiplication $\Delta$. For $l \in \frac{1}{2}\bbn$ and $ i,j \in I_l$, let
\begin{IEEEeqnarray*}{rCl}
f_j^{l} &=& {2l\choose l+j}_{|q|^2}^{1/2}\vartheta_R(y^{l-j}x^{l+j})={2l\choose l+j}_{|q|^2}^{1/2}a^{l-j}b^{l+j},\\
e_i^{l} &=& {2l\choose l+i}_{|q|^2}^{1/2}\vartheta_L(y^{l-i}x^{l+i})={2l \choose l+i}_{|q|^2}^{1/2}a^{l-i}c^{l+i},\\
V_l^R &=& \oplus_{j=-l}^l \bbc f_j^{(l)}, \qquad V_l^L=\oplus_{i=-l}^l \bbc e_i^{(l)},\\
T_l^R &=& \left.\Delta\right|_{V_l^R}, \qquad\qquad T_l^L=\left.\Delta\right|_{V_l^L}\,.
\end{IEEEeqnarray*}
Employing Propn. (\ref{rightaction},\ref{equi}) and the fact that $ \mathcal{O}(\bbc_q^2)_{2l}$ is an invariant subspace under $\psi_R$, one gets a corepresentation $T_l^R: V_l^R \longrightarrow V_l^R \otimes \cla_q$ of $\cla_q$ on $V_l^R$ which is equivalent to $T_l$. For $i,j\in I_l$, let $t_{ij}^l$ denote the matrix coefficients of $T_l^R$ with respect to the basis $\{f_j^l\}$ of $V_l^R$. Then, we have
\begin{IEEEeqnarray}{rCl}\label{eq2}
T_l^R(f_j^l)=\Delta(f_j^l)=\sum_{i=-l}^{l} f_i^l \otimes t_{ij}^l\,.
\end{IEEEeqnarray}
Similarly, let $w_{ij}^l$ be the matrix coefficients of $T_l^L$ with respect to the basis $\{e_i^l\}$ of $V_l^L$. Then, we have
\begin{IEEEeqnarray}{rCl} \label{eq(2)}
T_l^L(e_i^l)=\Delta(e_i^l)=\sum_{j=-l}^{l} w_{ij}^l \otimes  e_j^l \,.
\end{IEEEeqnarray}
Consider the sesquilinear forms $\langle \cdot, \cdot \rangle_R$ and $\langle \cdot, \cdot \rangle_L$ on $\cla_q$ defined as follows,
\begin{center}
$\langle x,y \rangle_L=h(x^*y)\quad,\quad\langle x,y \rangle_R=\overline{h(xy^*)} \qquad \qquad \mbox{ for } x,y \in \cla_q\,.$
\end{center}
Both the sesquilinear forms are positive definite, thanks to Thm. \ref{faithful}, and hence  $\cla_q$ is an inner product space under $\langle \cdot, \cdot \rangle_R$ and $\langle \cdot, \cdot \rangle_L\,$. The following proposition says that the decomposition of $\cla_q$ given in part (vi) of Propn. \ref{clafurther} is orthogonal. 
\bppsn\label{orthogonal}
Let $x \in \cla_q[m,n,r]$ and $y \in \cla_q[m^{'},n^{'},r^{'}]$. If $(m,n,r) \neq (m^{'},n^{'},r^{'})$ then,
\begin{center}
$\langle x,y \rangle_L=0  \qquad \mbox{ and } \qquad  \langle x,y \rangle_R=0\,.$
\end{center}
\eppsn 
\prf By Propn. \ref{clafurther}, we have $x^*\in\cla_q[-m,-n,-r]$ and $y^*\in\cla_q[-m^{'},-n^{'},-r^{'}]$. Hence, $x^*y \in\cla_q[m^{'}-m,n^{'}-n,r^{'}-r]$ and $xy^*\in\cla_q[m-m^{'},n-n^{'},r-r^{'}]$. Using Propn. \ref{clafurtherfurther}, it now follows that $h(x^*y)=h(xy^*)=0$.  
\qed
\medskip

\subsection{The matrix coefficients and Peter-Weyl decomposition}

For $x \in \cla_q$, define $\|x\|_L=\langle x,x\rangle_L^{1/2}$ and $\|x\|_R=\langle x,x\rangle_R^{1/2}$.
\blmma \label{induction}
Let $d_{r,s}=\|a^rb^s\|_L^2$. Then, one has the followings.
\begin{enumerate}
\item[(i)] $d_{r,s}=d_{r-1,s}-|q|^{2r}d_{r-1,s+1}\,\,\forall\,r\geq 1\,;$
\item[(ii)] $d_{r,s}=\Big(|q|^{(r+s)}|r+s+1|_{|q|}{r+s \choose s}_{|q|^2}\Big)^{-1}\,$.
\end{enumerate}
\elmma
\prf  For part $(i)$, we have
\begin{IEEEeqnarray*}{rCl}
d_{r,s} &=& h((b^*)^s(a^*)^ra^rb^s)=h((b^*)^sb^s(a^*)^ra^r)\\
&=& h\Big((b^*)^s(a^*)^{r-1}(1-|q|^2b^*b)a^{r-1}b^s\Big)\\
&=& h\Big((b^*)^s(a^*)^{r-1})a^{r-1}b^s\Big)-|q|^{2r}h\Big((b^*)^{s+1}(a^*)^{r-1}a^{r-1}b^{s+1}\Big)\\
&=& d_{r-1,s} -|q|^{2r}d_{r-1,s+1}\,.
\end{IEEEeqnarray*}
For part $(ii)$, we use induction based on the formula given in part $(i)$. If $r=0$, we have 
\begin{IEEEeqnarray*}{rCl}
 d_{0,s}&=& h((b^s)^*b^s)=\frac{1-|q|^2}{1-|q|^{2(s+1)}}=\frac{1}{1+|q|^2+|q|^4+ \cdots +|q|^{2s}}=\frac{1}{|q|^s|s+1|_{|q|}}\,.
\end{IEEEeqnarray*}
Hence, the claim is true for all the tuples $(0,s)$ where $s\in\bbn$. Assume that the claim holds for each tuple $(p,s)$, where $0\leq p<r$ and $s\in\bbn$. Using part $(i)$, we then have
\begin{IEEEeqnarray*}{rCl}
d_{r,s} &=& d_{r-1,s}-|q|^{2r}d_{r-1,s+1}\\
&=& \frac{1}{|q|^{r+s-1}|r+s|_{|q|} {r+s-1 \choose s}_{|q|^2}}-|q|^{2r} \frac{1}{|q|^{r+s}|r+s+1|_{|q|} {r+s \choose s+1}_{|q|^2}}\\
&=& \frac{(|q|^2,|q|^2)_{r-1}(|q|^2,|q|^2)_{s}}{(1+|q|^2+\cdots +|q|^{2r+2s-2})(|q|^2,|q|^2)_{r+s-1}}
-\frac{|q|^{2r}(|q|^2,|q|^2)_{r-1}(|q|^2,|q|^2)_{s+1}}{(1+|q|^2+\cdots +|q|^{2r+2s})(|q|^2,|q|^2)_{r+s}}\\
&=&\frac{(1-|q|^2)\prod_{i=1}^{r-1}(1-|q|^{2i})\prod_{i=1}^{s}(1-|q|^{2i})}{(1-q|^{2r+2s})\prod_{i=1}^{r+s-1}(1-|q|^{2i})} - \frac{|q|^{2r}(1-|q|^2)\prod_{i=1}^{r-1}(1-|q|^{2i})\prod_{i=1}^{s+1}(1-|q|^{2i})}{(1-q|^{2r+2s+2})\prod_{i=1}^{r+s}(1-|q|^{2i})}\\
&=&\frac{(1-|q|^2)\prod_{i=1}^{r-1}(1-|q|^{2i})\prod_{i=1}^{s}(1-|q|^{2i})}{\prod_{i=1}^{r+s}(1-|q|^{2i})}
\Big(1-\frac{|q|^{2r}(1-|q|^{2s+2})}{1-|q|^{2r+2s+2}}\Big)\\
&=&\frac{(1-|q|^2)\prod_{i=1}^{r}(1-|q|^{2i})\prod_{i=1}^{s}(1-|q|^{2i})}{(1-|q|^{2r+2s+2})\prod_{i=1}^{r+s}(1-|q|^{2i})}\\
&=& \frac{1}{(1+|q|^2+\cdots +|q|^{2r+2s}){r+s \choose s}_{|q|^2}}\\
&=& \frac{1}{|q|^{(r+s)}|r+s+1|_{|q|} {r+s \choose s}_{|q|^2}}\,.
\end{IEEEeqnarray*}
\qed 
 \blmma \label{induction1}
Let $c_{r,s}=\|a^rb^s\|_R^2$. Then, one has the followings.
\begin{enumerate}
\item[(i)] $c_{r,s}=\frac{1}{|q|^{2s}}c_{r-1,s}-\frac{1}{|q|^{2s}}c_{r-1,s+1}\,\,\forall\,r\geq 1\,;$
\item[(ii)] $c_{r,s}=|q|^{2r}\Big(|q|^{(r+s)}|r+s+1|_{|q|} {r+s \choose s}_{|q|^2}\Big)^{-1}=|q|^{2r}d_{r,s}\,$.
\end{enumerate}
\elmma
\prf  For part $(i)$, we have
\begin{IEEEeqnarray*}{rCl}
c_{r,s} &=& \overline{h(a^rb^s(b^*)^s(a^*)^r)}\\
&=& \frac{1}{|q|^{2s}}\overline{h(a^{r-1}b^s(b^*)^saa^*(a^*)^{r-1})}\\
&=& \frac{1}{|q|^{2s}}\overline{h(a^{r-1}b^s(b^*)^s(a^*)^{r-1})}-\frac{1}{|q|^{2s}}\overline{h(a^{r-1}b^s(b^*)^sbb^*(a^*)^{r-1})}\\
&=& \frac{1}{|q|^{2s}}c_{r-1,s}-\frac{1}{|q|^{2s}}c_{r-1,s+1}\,.
\end{IEEEeqnarray*}
For part $(ii)$, we use induction based on the formula given in part $(i)$. If $r=0$, we have 
\begin{IEEEeqnarray*}{rCl}
c_{0,s}&=& \overline{h(b^s(b^*)^s)}=\Big(|q|^s|s+1|_{|q|}\Big)^{-1}\,.
\end{IEEEeqnarray*}
Observe that this is same as $d_{0,s}$ in the previous lemma due to normality of $b$. Hence, the claim is true for all the tuples $(0,s)$ where $s\in\bbn$. Assume that the claim holds for each tuple $(p,s)$, where $0\leq p<r$ and $s\in\bbn$. Then, using part $(i)$, we have
\begin{IEEEeqnarray*}{rCl}
c_{r,s}&=&\frac{1}{|q|^{2s}}c_{r-1,s}-\frac{1}{|q|^{2s}}c_{r-1,s+1}\\
&=& \frac{|q|^{2r-2}}{|q|^{2s}|q|^{(r+s-1)}|r+s|_{|q|} {r+s-1 \choose s}_{|q|^2}}- \frac{|q|^{2r-2}}{|q|^{2s}|q|^{(r+s)}|r+s+1|_{|q|} {r+s \choose s+1}_{|q|^2}}\\
&=& \frac{1-|q|^2}{|q|^{2s}}\Big(\frac{\prod_{l=1}^{r-1}(1-|q|^{2l})|q|^{2(r-1)}}{\prod_{l=1}^r(1-|q|^{2(l+s)})}
  -\frac{\prod_{l=1}^{r-1}(1-|q|^{2l})|q|^{2(r-1)}}{\prod_{l=1}^r(1-|q|^{2(l+s+1)})}\Big)\\
&=& \frac{(1-|q|^2)|q|^{2(r-1)}}{|q|^{2s}}\frac{\prod_{l=1}^{r-1}(1-|q|^{2l})(|q|^{2(s+1)}-|q|^{2(r+s+1)})}{\prod_{l=1}^{r+1}(1-|q|^{2(l+s)})}\\
&=& \frac{(1-|q|^2)|q|^2\prod_{l=1}^r(1-|q|^{2l})|q|^{2(r-1)}}{\prod_{l=1}^{r+1}(1-|q|^{2(l+s)})}\\
&=& \frac{(1-|q|^2)\prod_{l=1}^r(1-|q|^{2l})|q|^{2r}}{\prod_{l=1}^{r+1}(1-|q|^{2(l+s)})}\\
&=& \frac{|q|^{2r}}{|q|^{(r+s)}|r+s+1|_{|q|}{r+s\choose s}_{|q|^2}}\,.
\end{IEEEeqnarray*}
Thus, $c_{r,s}=|q|^{2r}d_{r,s}$ by Lemma \ref{induction}.\qed 
\bppsn\label{orthonormal basis}
In the  Hilbert space $(V_l^R, \langle \cdot, \cdot \rangle_L)$, one has the followings.
\begin{enumerate}[(i)]
 \item $\|f_j^l\|=\frac{1}{|q|^{l}|2l+1|_{|q|}^{1/2}}\,;$
 \item The ordered set $B_l:=\big\{|q|^{l}|2l+1|_{|q|}^{1/2}f_j^l\big\}_{i \in I_l}$ is an orthonormal basis of $V_l^R\,;$
 \item The matrix coefficients of $T_l^R$ with respect to $B_l$ are $t_{ij}^l, i,j \in I_l\,.$
\end{enumerate}
\eppsn
\prf The first part follows from Lemma \ref{induction}. For $r\neq r^{'}$, using the faithful representations defined in \ref{representation}, we have for $n\geq r$
\begin{center}
$\langle v_{n,0,0}\,,\,a^{r^{'}}b^{2l-r^{'}}(b^*)^{2l-r}(a^*)^r(v_{n,0,0})\rangle=\langle v_{n,0,0}\,,\,Cv_{n-r+r^{'},r-r^{'},0}\rangle=0$
\end{center}
for orthonormal basis $\{v_{n,m,k}\}$ of $\ell^2(\mathbb{N})\otimes\ell^2(\mathbb{Z})\otimes\ell^2(\mathbb{Z})$, where $C$ is a nonzero constant. Since the above inner-product is automatically zero for $0\leq n<r$, by Thm. \ref{haar} we get that
\begin{center}
$\langle a^{r^{'}}b^{2l-r^{'}},\,a^rb^{2l-r}\rangle_R=\overline{h\Big(a^{r^{'}}b^{2l-r^{'}}(b^*)^{2l-r}(a^*)^r\Big)}=0\,.$
\end{center}
Combining this with part $(i)$, we get the claim. The last part follows from eqn. \ref{eq2}.\qed
\blmma For $m \geq 0$, one has the following~:
\begin{IEEEeqnarray*}{rCl}
 a^m(a^*)^m &=& \sum_{k= 0}^m (-1)^k {m \choose k}_{|q|^2} |q|^{k^2+k-2km}(bb^*)^k\,,\\
 (a^*)^ma^m &=& \sum_{k= 0}^m (-1)^k {m \choose k}_{|q|^2} |q|^{k^2+k}(bb^*)^k\,.
\end{IEEEeqnarray*}
\elmma
\prf Follows from induction on $m$.\qed
\bppsn\label{matrixcoefficients}
For $l \in \frac{1}{2}\bbn$, let   $t_{ij}^l$, $i,j \in I_l$ be the matrix coefficients of the corepresentation $T_l^R$ with 
respect to the basis $\{f_j^l: j \in I_l\}$. Then, we have
\begin{enumerate}
\item[(i)]
\begin{IEEEeqnarray*}{rCl}
t_{ij}^l= \sum_{\substack{m+n=l-i\\ 0\leq m \leq l-j\\0\leq n \leq l+j}}q^{n(l-j-m)}
\frac{{ 2l \choose l+j}_{|q|^2}^{1/2}}{{ 2l \choose l+i}_{|q|^2}^{1/2}}{l-j \choose m}_{|q|^2}{l+j \choose n}_{|q|^2}a^mc^{l-j-m}b^nd^{l+j-n}\,.
\end{IEEEeqnarray*}
\item[(ii)] For 
 $m,n \in \bbn$ satisfying $  0\leq m \leq l-j, 0\leq n \leq l+j$ and $m+n=l-i$, let 
 \[
   C_q(l,i,j,m,n)=  \frac{|q|^{(l-j-m)(l-j-m+2n+1)+2n(l+i)}}{(\bar{q})^{(i-j)(l+i)}}
 \frac{{ 2l \choose l+j}_{|q|^2}^{1/2}}{{ 2l \choose l+i}_{|q|^2}^{1/2} } {l-j \choose m}_{|q|^2}{l+j \choose n}_{|q|^2}.
 \]
Then, one has 
 \[
 t_{ij}^l= \sum_{\substack{m+n=l-i\\ 0\leq m \leq l-j\\0\leq n \leq l+j}}(-1)^{l-j-m}C_q(l,i,j,m,n)a^m(a^*)^{l+j-n}(b^*)^{l-j-m}b^nD^{l+i}\,.
 \]
\item[(iii)] For $\zeta=bb^*$ and $e_{m,n}$ as defined in $\ref{defn of e_m,n}$, one has
\[
t_{i,j}^l=e_{-2i,-2j}P_{i,j}(\zeta)D^{l+i}
\]
for some polynomial $P_{i,j}$ where, $deg(P_{i,j})=l-\max \{|i|,|j|\}$.
\end{enumerate}
\eppsn
\prf For $c=-\bar{q}Db^*$ and $d=Da^*$, we have
\begin{IEEEeqnarray*}{rCl}
\Delta(a)=a \otimes a+b \otimes c\,,\qquad\Delta(b)=a\otimes b+b \otimes d\,. 
\end{IEEEeqnarray*}
Using Propn. $(2)$ in (page $39$, eqn. $18$ in \cite{KliSch-1997aa}), we have 
\begin{IEEEeqnarray*}{rCl}
\Delta(a^{l-j})&= &\sum_{m=0}^{l-j}{l-j \choose m}_{|q|^2} a^mb^{l-j-m}\otimes a^mc^{l-j-m}\\
\Delta(b^{l+j})&=& \sum_{n=0}^{l+j}{l+j \choose n}_{|q|^2} a^nb^{l+j-n}\otimes b^nd^{l+j-n}\\
\end{IEEEeqnarray*}
since, $ab=q^{-1}ba,\,ac=(\bar{q})^{-1}ca,\,bd=(\bar{q})^{-1}db\,$. Multiplying these two expressions, we get that
\begin{IEEEeqnarray*}{rCl}
\Delta(f_j^l)&=&\Delta\Big({ 2l \choose l+j}_{|q|^2}^{1/2}  a^{l-j}b^{l+j}\Big)\\
&=& \sum_{m=0}^{l-j}\sum_{n=0}^{l+j}{ 2l \choose l+j}_{|q|^2}^{1/2} {l-j \choose m}_{|q|^2}{l+j \choose n}_{|q|^2}a^mb^{l-j-m}a^nb^{l+j-n}\otimes a^mc^{l-j-m}b^nd^{l+j-n}\\
&=& \sum_{m=0}^{l-j}\sum_{n=0}^{l+j}q^{n(l-j-m)}{ 2l \choose l+j}_{|q|^2}^{1/2} {l-j \choose m}_{|q|^2}{l+j \choose n}_{|q|^2}a^{m+n}b^{2l-m-n}\otimes a^mc^{l-j-m}b^nd^{l+j-n}\,.
\end{IEEEeqnarray*}
Putting $\,m+n=l-i$, we get that
\begin{IEEEeqnarray*}{rCl}
\Delta(f_j^l)&=& \sum_{\substack{m+n=l-i\\ 0\leq m\leq l-j\\0\leq n\leq l+j}}q^{n(l-j-m)}{ 2l \choose l+j}_{|q|^2}^{1/2} {l-j \choose m}_{|q|^2}{l+j \choose n}_{|q|^2}a^{l-i}b^{l+i}\otimes a^mc^{l-j-m}b^nd^{l+j-n}\\
&=&\sum_{\substack{m+n=l-i\\ 0\leq m \leq l-j\\0\leq n \leq l+j}}q^{n(l-j-m)}\frac{{ 2l \choose l+j}_{|q|^2}^{1/2}}{{ 2l \choose l+i}_{|q|^2}^{1/2} } {l-j \choose m}_{|q|^2}{l+j \choose n}_{|q|^2}f_i^l\otimes a^mc^{l-j-m}b^nd^{l+j-n}\,.
\end{IEEEeqnarray*}
Comparing with eqn. \ref{eq2}, we get that
\begin{IEEEeqnarray*}{lCl}
t_{ij}^l &=& \sum_{\substack{m+n=l-i\\ 0\leq m \leq l-j\\0\leq n \leq l+j}}q^{n(l-j-m)}\frac{{ 2l \choose l+j}_{|q|^2}^{1/2}}{{ 2l \choose l+i}_{|q|^2}^{1/2} } {l-j \choose m}_{|q|^2}{l+j \choose n}_{|q|^2}a^mc^{l-j-m}b^nd^{l+j-n}\\
&=& \sum_{\substack{m+n=l-i\\ 0\leq m \leq l-j\\0\leq n \leq l+j}}q^{n(l-j-m)}\frac{{ 2l \choose l+j}_{|q|^2}^{1/2}}{{ 2l \choose l+i}_{|q|^2}^{1/2}}{l-j \choose m}_{|q|^2}{l+j \choose n}_{|q|^2}a^m(-\bar{q}Db^*)^{l-j-m}b^n(Da^*)^{l+j-n}\\
&=& \sum_{\substack{m+n=l-i\\ 0\leq m \leq l-j\\0\leq n \leq l+j}}(-1)^{l-j-m}C_q(l,i,j,m,n)a^m(a^*)^{l+j-n}(b^*)^{l-j-m}b^nD^{l+i}\,.
\end{IEEEeqnarray*}
This completes part $(i)$ and part $(ii)$. To get the last part, one needs to consider four cases namely, $(i)\,i+j\geq 0,\,i \geq j\quad (ii)\,i+j\geq 0,\,i\leq j\quad (iii)\,i+j\leq 0,\,i \geq j$ and $(iv)\,i+j\leq 0,\,i\leq j$. We will prove the claim for the first case only as all the other cases are similar. 
\begin{IEEEeqnarray*}{lCl}
t_{ij}^l &=& \sum_{\substack{m+n=l-i\\ 0\leq m \leq l-i\\0\leq n \leq l-i}}(-1)^{l-j-m}C_q(l,i,j,m,n)a^m(a^*)^{l+j-(l-i-m)}(b^*)^{l-j-m}b^{l-i-m}D^{l+i}\\
&=& \sum_{\substack{m+n=l-i\\ 0\leq m \leq l-i\\0\leq n \leq l-i}}(-1)^{l-j-m}C_q(l,i,j,m,n) a^m(a^*)^m(a^*)^{i+j}(b^*)^{i-j}(b^*)^{l-i-m}b^{l-i-m}D^{l+i}\\
&=& \sum_{\substack{m+n=l-i\\ 0\leq m \leq l-i\\0\leq n \leq l-i}}(-1)^{l-j-m}(\overline{q})^{(i+j)(i-j)}C_q(l,i,j,m,n) a^m(a^*)^m\zeta^{l-i-m}e_{-2i,-2j}D^{l+i}\\
&=& \Big(\sum_{\substack{m+n=l-i\\ 0\leq m\leq l-i\\0\leq n \leq l-i}}\sum_{k= 0}^m(-1)^{l-j-m+k}(\overline{q})^{(i+j)(i-j)}|q|^{2(i+j)(l-i-m)}C_q(l,i,j,m,n){m \choose k}_{|q|^2}\\
& & \qquad\qquad\qquad |q|^{k^2+k-2km}\zeta^{l-i-m+k}\Big)e_{-2i,-2j}D^{l+i}\,.
\end{IEEEeqnarray*}
Hence, $t_{ij}^l=R_{i,j}(\zeta)e_{-2i,-2j}D^{l+i}$ where $R_{i,j}$ is a polynomial of degree $l-i$. Applying Propn. \ref{deg same} we get the claim.\qed
\brmrk
A more concrete description of $t^l_{ij}$ in terms of the little $q$-Jacobi polynomials is given in the next subsection $4.3$.
\ermrk
\blmma \label{class}
For $l\in\frac{1}{2}\bbn,$ one has the followings.
\begin{enumerate}[(i)]
\item $f_j^l = t_{-l,j}^l$ and $e_i^l = t_{i,-l}^l$  for $i \in I_l\,;$
\item $t_{ij}^l =w_{ij}^l $ for  $i,j \in I_l\,;$
\item $t_{ij}^l\in\mathcal{A}_q[-2i,-2j, l+i]$  for $i,j \in I_l\,.$
\end{enumerate}
\elmma 
\prf Part $(i)$ and $(iii)$ are immediate consequences of Propn. \ref{matrixcoefficients}. To show part $(ii)$, observe that 
\begin{center}
$\Delta(t_{ij}^l)=\sum_{k \in I_l}t_{ik}^l \otimes t_{kj}^l\,.$
\end{center}
Taking $j=-l$ and using the fact that $e_i^l = t_{i,-l}^l$, we get that
\begin{center}
$\Delta(e_i^l)=\sum_{k \in I_l}t_{ik}^l \otimes e_k^l\,.$
\end{center}
Comparing with the eqn. \ref{eq(2)}, we get the claim.\qed
\blmma \label{leftright}
Let $m,n \in \bbz$ such that $m-n$ is even. Then, for $\,k_1,k_2\in\bbn\,,$
\begin{enumerate}[(i)]
\item $\langle e_{m,n}\zeta^{k_1},e_{m,n}\zeta^{k_2}\rangle_L=|q|^{-m-n} \langle e_{m,n}\zeta^{k_1},e_{m,n}\zeta^{k_2}\rangle_R\,\,,$ where $e_{m,n}$ is as defined in $\ref{defn of e_m,n}\,.$
\item For $x \in \cla_q[m,n]$, one has $\langle x,x \rangle_L= |q|^{-m-n} \langle x,x \rangle_R\,.$
\end{enumerate}
\elmma
\prf Let $m+n\geq 0$ and $m\geq n$. Using the expressions for $d_{r,s}$ and $c_{r,s}$ from the Lemmas (\ref{induction},\ref{induction1}) and the fact that $ba^r(a^*)^r=a^r(a^r)^*b$, we have 
\begin{IEEEeqnarray*}{rCl}
 \langle e_{m,n}\zeta^{k_1},e_{m,n}\zeta^{k_2}\rangle_L&=& 
 \big\langle a^{\frac{m+n}{2}}b^{\frac{m-n}{2}}\zeta^{k_1},a^{\frac{m+n}{2}}b^{\frac{m-n}{2}}\zeta^{k_2}\big\rangle_L\\
 &=& h\big(\zeta^{k_1}(b^*)^{\frac{m-n}{2}}(a^*)^{\frac{m+n}{2}}a^{\frac{m+n}{2}}b^{\frac{m-n}{2}}\zeta^{k_2}\big)\\
 &=& h\big((b^*)^{\frac{m-n}{2}+k_1+k_2}(a^*)^{\frac{m+n}{2}}a^{\frac{m+n}{2}}b^{\frac{m-n}{2}+k_1+k_2}\big)\\
 &=& d_{\frac{m+n}{2},\frac{m-n}{2}+k_1+k_2}\\
 &=& |q|^{-m-n}\,c_{\frac{m+n}{2},\frac{m-n}{2}+k_1+k_2}\\
 &=& |q|^{-m-n}\,\overline{h\big(a^{\frac{m+n}{2}}b^{\frac{m-n}{2}+k_1+k_2}(b^*)^{\frac{m-n}{2}+k_1+k_2}(a^*)^{\frac{m+n}{2}}\big)}\\
 &=& |q|^{-m-n}\langle e_{m,n}\zeta^{k_1},e_{m,n}\zeta^{k_2}\rangle_R\,.
\end{IEEEeqnarray*}
In the case of $m+n\leq 0$ and $m\leq n$, we have
\begin{IEEEeqnarray*}{rCl}
 \langle e_{m,n}\zeta^{k_1},e_{m,n}\zeta^{k_2}\rangle_L&=&
 \big\langle (b^*)^{\frac{n-m}{2}}(a^*)^{\frac{-m-n}{2}}\zeta^{k_1},(b^*)^{\frac{n-m}{2}}(a^*)^{\frac{-m-n}{2}}\zeta^{k_2}\big\rangle_L\\
 &=& h\big(\zeta^{k_1}a^{\frac{-m-n}{2}}b^{\frac{n-m}{2}}(b^*)^{\frac{n-m}{2}}(a^*)^{\frac{-m-n}{2}}\zeta^{k_2}\big)\\
 &=& |q|^{k_1(-m-n)}|q|^{k_2(-m-n)}h\big( a^{\frac{-m-n}{2}}b^{\frac{n-m}{2}+k_1+k_2}(b^*)^{\frac{n-m}{2}+k_1+k_2}(a^*)^{\frac{-m-n}{2}}\big)\\
 &=&|q|^{k_1(-m-n)}|q|^{k_2(-m-n)}\overline{c_{\frac{-m-n}{2},\frac{n-m}{2}+k_1+k_2}}\\
 &=& |q|^{-m-n}|q|^{k_1(-m-n)}|q|^{k_2(-m-n)}\overline{d_{\frac{-m-n}{2},\frac{n-m}{2}+k_1+k_2}}\\
 &=& |q|^{-m-n}|q|^{k_1(-m-n)}|q|^{k_2(-m-n)}\overline{h\big((b^*)^{\frac{n-m}{2}+k_1+k_2}(a^*)^{\frac{-m-n}{2}}a^{\frac{-m-n}{2}}b^{\frac{n-m}{2}+k_1+k_2}\big)}\\
 &=&|q|^{-m-n}\,\overline{h\big((b^*)^{\frac{n-m}{2}}(a^*)^{\frac{-m-n}{2}} \zeta^{k_1+k_2}a^{\frac{-m-n}{2}}b^{\frac{n-m}{2}}\big)}\\
&=&|q|^{-m-n}   \langle e_{m,n}\zeta^{k_1},e_{m,n}\zeta^{k_2}\rangle_R\,. 
\end{IEEEeqnarray*}
The remaining two cases follow along the line of the above calculations. This proves the first part of the claim. To prove the second part, take a polynomial $P(\zeta)=\sum_{k=0}^n\alpha_k\zeta^k$. Invoking the first part of this Lemma, we have  
\begin{IEEEeqnarray*}{rCl}
\langle e_{m,n}P(\zeta),e_{m,n}P(\zeta)\rangle_L
&=&\sum_{k_1=0}^n\sum_{k_2=0}^n \overline{\alpha_{k_1}}\alpha_{k_2}  \langle e_{m,n} \zeta^{k_1}, e_{m,n}\zeta^{k_2}\rangle_L\\
&=&|q|^{-m-n} \langle e_{m,n}P(\zeta),e_{m,n}P(\zeta)\rangle_R\,.
\end{IEEEeqnarray*}
For $r\in\bbz$, from the commutation relations in \ref{relations}, we have 
\begin{center}
$(D^r)^*\zeta^{k_1}e_{m,n}^*e_{m,n}\zeta^{k_1}D^r= \zeta^{k_1}e_{m,n}^*e_{m,n}\zeta^{k_1}\,.$
\end{center}
Hence, we get that
\begin{IEEEeqnarray}{rCl} \label{e1}
  \langle e_{m,n}P(\zeta)D^r,e_{m,n}P(\zeta)D^r\rangle_L
  &=&|q|^{-m-n} \langle e_{m,n}P(\zeta)D^r,P(\zeta)e_{m,n}D^r\rangle_R\,.
\end{IEEEeqnarray}
If $r\neq s$, we have from Propn. \ref{clafurther} that $e_{m,n}P(\zeta)D^r(e_{m,n}P(\zeta)D^s)^* \in \cla_q[0,0,r-s]$ and $(e_{m,n}P(\zeta)D^r)^*e_{m,n}P(\zeta)D^s \in \cla_q[0,0,s-r]$. Using Propn. $\ref{clafurtherfurther}$, we get that
\begin{IEEEeqnarray}{rCl} \label{e2}
\langle e_{m,n}P(\zeta)D^r,e_{m,n}P(\zeta)D^s\rangle_L=0= \langle e_{m,n}P(\zeta)D^r,e_{m,n}P(\zeta)D^s\rangle_R\,.
\end{IEEEeqnarray}
Take $x\in \cla_q[m,n]$. By part $(v)$ of Propn. $\ref{clafurther}$, we have $x=\sum_{r=-k}^ke_{m,n}P_r(\zeta)D^r$ for some polynomials $P_r\,,\,-k \leq r\leq k,\,k\in\bbn$. Using eqns. $(\ref{e1},\ref{e2})$, we finally have
\begin{IEEEeqnarray*}{rCl}
\langle x,x\rangle_L&=& \sum_{r=-k}^k  \langle e_{m,n}P_r(\zeta)D^r,e_{m,n}P_r(\zeta)D^r\rangle_L\\
&=&|q|^{-m-n}\sum_{r=-k}^k  \langle e_{m,n}P_r(\zeta)D^r,e_{m,n}P_r(\zeta)D^r\rangle_R\\
&=& |q|^{-m-n}\langle x,x\rangle_R\,.
\end{IEEEeqnarray*}
\qed
\brmrk
Using the relations in \ref{relations}, it is easy to check that the map $\sigma:\cla_q\longrightarrow\cla_q$ defined by $\,a\rightarrow q^{-2}a,\,a^*\rightarrow q^2a^*,\,b\rightarrow b,\,b^*\rightarrow b^*,\,D\rightarrow D,\,D^*\rightarrow D^*$ extends to an algebra automorphism. One can now follow (Lemma $3.3$ in \cite{Masuda-1991aa}, or \cite{KliSch-1997aa}) and show that $h(xy)=h(\sigma(y)x)$ for $x,y\in\cla_q,\,$i,e. $\sigma$ is the modular automorphism. Using this, the above lemma can be proved. However, the proof given here is more straightforward and uses only norms given in Lemmas (\ref{induction},\ref{induction1}). 
\ermrk

\bthm \label{irreduciblerep}
For $l\in\frac{1}{2}\bbn,\,T_l^R$ is an irreducible corepresentation of $\cla_q\,$. Moreover, the matrix $\big(\big(t_{ij}^l\big)\big)$ of $T_l^R$ with respect to the orthonormal basis $B_l$ in Propn. $\ref{orthonormal basis}$ is a unitary element in $M_{2l+1}(\bbc)\otimes \cla_q\,$. 
\ethm
\prf Tha claim of unitarizability of $T_l^R$ with respect to the inner-product $\langle \cdot,  \cdot \rangle_L$ follows along the same line of argument in (Propn. $3.5$, \cite{Masuda-1991aa}). Since $B_l$ is an orthonormal basis of the inner-product space $(V_l^R,\langle \cdot,  \cdot \rangle_L)$, the matrix of $T_l^R$ with respect to $B_l$ is a unitary matrix. Irreducibility of $T_l^R$  follows along the same line of argument in the proof of Thm. $9$ in (\cite{KliSch-1997aa}, page $110$).\qed
\brmrk One can also use Propn. \ref{matrixcoefficients} and directly verify that $S\big(t_{ij}^l\big)=\big(t_{ji}^l\big)^*$, in order to show that $\big(t_{ij}^l\big)$ is a unitary element in $M_{2l+1}(\bbc)\otimes \cla_q\,$.
\ermrk
For $m\in\bbz$, let $V_l^RD^m=\bigoplus_{i=-l}^l\bbc f_i^{(l)}D^m$ and $B_{l,m}=\{|q|^{l}|2l+1|_{|q|}f_i^{(l)}D^m: i \in I_l\}$.  
Define a representation $T_l^RD^m:V_l^RD^m\longrightarrow V_l^RD^m \otimes \cla_q$ by a linear map 
\[
 T_l^RD^m(f_i^lD^m)= \sum_{j\in I_l} f_j^lD^m \otimes t_{ji}^lD^m.
\]
\bthm For $m \in \bbz$ and $l \in \frac{1}{2}\bbn$, 
$ T_l^RD^m$ is an irreducible  corepresentation of $\cla_q$ and its
matrix coefficients with respect to the basis $B^l$ are $\{t_{ij}^lD^m:i,j \in I_l\}$. Moreover, the matrix $\big(\big(t_{ij}^lD^m\big)\big)$ is a unitary element of $M_{2l+1}(\bbc)\otimes \cla_q\,$. 
\ethm
\prf For $i,j \in I_l$, we have 
\[
 \Delta(t_{ij}^lD^m)=\Delta(t_{ij}^l)\Delta(D^m)=(\sum_{k \in I_l}t_{ik}^l \otimes t_{kj}^l)(D^m \otimes D^m)
 =\sum_{k \in I_l}t_{ik}^lD^m \otimes t_{kj}^lD^m.
\]
This proves that $T_l^RD^m$ is a corepresentation of $\cla_q$. Since  $T_l^RD^m(f_i^lD^m)= \sum_{j\in I_l} f_j^lD^m \otimes t_{ji}^lD^m$, it follows that 
matrix coefficients with respect to the basis $B_{l,m}$ are $\{t_{ij}^lD^m:i,j \in I_l\}$. Moreover, we have
\[
 S(t_{ij}^lD^m)=S(D^m)S(t_{ij}^l)=(D^m)^*(t_{ji}^l)^*=(t_{ji}^lD^m)^*
\]
which shows that $T_l^RD^m$ is a unitary corepresentation. Irreducibility follows by a similar argument used to prove  Thm. $9$ in \cite{KliSch-1997aa} (see page no. 110, \cite{KliSch-1997aa}). 
\qed 

\bthm\label{Peter}
\textbf{Peter-Weyl decomposition~:} Given a corepresentation $T$ of $\cla_q$, let $\mathcal{C}(T)$ denotes the 
vector subspace of $\cla_q$ generated by the matrix coefficients of $T$. Then, we have
\begin{enumerate}[(i)]
\item $\cla_q=\bigoplus_{l \in \frac{1}{2}\bbn, m \in \bbz}\quad \mathcal{C}(T_l^RD^m)$.
\item 
\begin{IEEEeqnarray}{rCl}
\langle t_{ij}^{l}D^m, t_{i^{'}j^{'}}^{l^{'}}D^{m^{'}}\rangle_R&=&|q|^{-2j}|2l+1|_{|q|}^{-1}\delta_{ll^{'}}\delta_{ii^{'}}\delta_{jj^{'}}\delta_{mm^{'}}\\
\langle t_{ij}^{l}D^m, t_{i^{'}j^{'}}^{l^{'}}D^{m^{'}}\rangle_L&=&|q|^{2i}|2l+1|_{|q|}^{-1}\delta_{ll^{'}}\delta_{ii^{'}}\delta_{jj^{'}}\delta_{mm^{'}}
\end{IEEEeqnarray}
\item The set $\{T_l^RD^m: l \in \frac{1}{2}\bbn, m \in \bbz\}$ is a complete list of irreducible mutually inequivalent corepresentations of $\cla_q\,$.
\item The set $\big\{|q|^{-i}|2l+1|_{|q|}^{1/2}t_{ij}^{l}D^m:l \in \frac{1}{2}\bbn, m \in \bbz\big\}$ is an orthonormal basis of $L^2(h)$.
\end{enumerate}
\ethm
\prf 
\textbf{proof of part (i):} Fix $i,j \in \frac{1}{2}\bbz$. To prove the claim, it suffices to show that
\begin{center}
$\cla_q[-2i,-2j]\subseteq\bigoplus_{l \in \frac{1}{2}\bbn, m \in \bbz}\quad \mathcal{C}(T_l^RD^m)\,.$
\end{center}
Since, $ \mathcal{C}(T_l^RD^m)= \mathcal{C}(T_l^R)D^m$ and $\cla_q[-2i,-2j,r]=\cla_q[-2i,-2j,0]D^r$ for $r \in \bbz$, 
it is enough to show that
\begin{center}
$\cla_q[-2i,-2j,0]\subseteq\bigoplus_{l\in\frac{1}{2}\bbn, m\in\bbz}\quad \mathcal{C}(T_l^RD^m)\,.$
\end{center}
By Proposition $\ref{clafurther}$, we have $\cla_q[-2i,-2j,0]= e_{-2i,-2j}\bbc[\zeta]$. Therefore, given any $n\in\bbn$ if we 
can exhibit an element in $\bigoplus_{l\in\frac{1}{2}\bbn, m\in\bbz}\,\,\mathcal{C}(T_l^RD^m)$ of the form 
$e_{-2i,-2j}P(\zeta)$, where $P$ is a polynomial of degree  $n$, then we get the claim. For that, choose $l_n=n+\max\{|i|,|j|\}$ and $m_n= -l_n-i$. By part $(ii)$ of Propn. \ref{matrixcoefficients}, we have
\begin{center}
$t_{i,j}^{l_n}D^{m_n}=e_{-2i,-2j}P_{i,j}(\zeta)$
\end{center}
where, $P_{i,j}$ is a polynomial of degree $n$ and this completes the proof.\\
\textbf{proof of part (ii):} With Lemma \ref{leftright} at our disposal, we can follow the proof of Thm. $17$ in (\cite{KliSch-1997aa}, page $115$).\\
\textbf{proof of part (iii):} From part $(i)$ of the claim, it follows that $\{T_l^RD^m:l\in\frac{1}{2}\bbn, m\in\bbz\}$ is a complete list of irreducible corepresentations of $\cla_q\,$. Also, if two corepresentations $T_1$ and $T_2$ are equivalent then $\mathcal{C}(T_1)=\mathcal{C}(T_1)$. By part $(ii)$, 
we have $\mathcal{C}(T_l^RD^m)\perp\mathcal{C}(T_{l^{'}}^RD^{m^{'}})$ if $(l,m)\neq (l^{'},m^{'})$. This proves the assertion.
\\
\textbf{proof of part (iv):} Follows from part $(i)$ and $(ii)$. 
\qed

\subsection{The little \texorpdfstring{$q$-}{}Jacobi polynomials and matrix coefficients}

Recall that the little $q$-Jacobi polynomials are defined as the following,
\begin{IEEEeqnarray*}{rCl}
\mathcal{P}_n^{(\alpha,\beta)}(z;q) & = & \sum_{r\geq 0}\frac{(q^{-n};q)_r(q^{\alpha+\beta+n+1};q)_r}{(q;q)_r(q^{\alpha+1};q)_r}(qz)^r\,,
\end{IEEEeqnarray*}
where, $(a;q)_m:=\prod_{k=0}^{m-1}(1-aq^k)$ for any $a\in\mathbb{C}$.
\bthm
The matrix coefficients $t_{ij}^lD^k$ are expressed in terms of the little q-Jacobi polynomials in the following way~:
\begin{enumerate}[(i)]
\item for the case of $\,i+j\leq 0,\,i\geq j$,
\begin{center}
$a^{-(i+j)}c^{i-j}(\bar{q})^{(j-i)(l+j)}\frac{{2l\choose l+j}^{1/2}_{|q|^2}}{{2l\choose l+i}^{1/2}_{|q|^2}}{l-j\choose i-j}_{|q|^2}\mathcal{P}_{l+j}^{(i-j,-i-j)}(bb^*;|q|^2)D^{l+j+k}\,;$
\end{center}
\item for the case of $\,i+j\leq 0,\,i\leq j$,
\begin{center}
$a^{-(i+j)}b^{j-i}q^{(i-j)(l+i)}\frac{{2l\choose l+j}^{1/2}_{|q|^2}}{{2l\choose l+i}^{1/2}_{|q|^2}}{l+j\choose j-i}_{|q|^2}\mathcal{P}_{l+i}^{(j-i,-i-j)}(bb^*;|q|^2)D^{l+i+k}\,;$
\end{center}
\item for the case of $\,i+j\geq 0,\,i\leq j$,
\begin{center}
$q^{(i-j)(l+i)}\frac{{2l\choose l+j}^{1/2}_{|q|^2}}{{2l\choose l+i}^{1/2}_{|q|^2}}{l+j\choose j-i}_{|q|^2}\mathcal{P}_{l-j}^{(j-i,i+j)}(bb^*;|q|^2)(a^*)^{i+j}b^{j-i}D^{l+i+k}\,;$
\end{center}
\item for the case of $\,i+j\geq 0,\,i\geq j$,
\begin{center}
$(\bar{q})^{(j-i)(l+j)}\frac{{2l\choose l+j}^{1/2}_{|q|^2}}{{2l\choose l+i}^{1/2}_{|q|^2}}{l-j\choose i-j}_{|q|^2}\mathcal{P}_{l-i}^{(i-j,i+j)}(bb^*;|q|^2)(a^*)^{i+j}c^{i-j}D^{l+j+k}\,.$
\end{center}
\end{enumerate}
\ethm
\prf The proof is purely computational and we only mention the key steps of the computations for Case $(i)$. All the other cases follow from similar computations. First, by induction one verifies that $(Db^*)^n=\Big(\frac{\bar{q}}{|q|}\Big)^{n^2-n}D^n(b^*)^n$. Then, using this and the defining relations in \ref{relations}, the following holds,
\begin{IEEEeqnarray*}{rCl}
&  & a^mc^{l-j-m}b^nd^{l+j-n}\\
& = & a^mc^{i-j+n}b^nd^{l+j-n}\\
& = &  (-\bar{q})^{i-j+n}\Big(\frac{\bar{q}}{|q|}\Big)^{n^2-n}a^m(Db^*)^{i-j}D^n(b^*)^nb^n(a^*)^{l+j-n}D^{l+j-n}\\
& = & (-\bar{q})^n\Big(\frac{\bar{q}}{|q|}\Big)^{n^2-n}a^mc^{i-j}D^n(bb^*)^n(a^*)^{i+j+m}D^{l+j-n}\\
& = & (-\bar{q})^n\Big(\frac{\bar{q}}{|q|}\Big)^{n^2-n}(|q|)^{-2n(i+j+m)}a^mc^{i-j}(a^*)^{i+j+m}D^{l+j}(bb^*)^n\\
& = & (-\bar{q})^n\Big(\frac{\bar{q}}{|q|}\Big)^{n^2-n}(|q|)^{-2n(i+j+m)}(\bar{q})^{-(i-j)(i+j+m)}a^{-(i+j)}c^{i-j}a^{i+j+m}(a^*)^{i+j+m}D^{l+j}(bb^*)^n\\
& = & (-\bar{q})^n\Big(\frac{\bar{q}}{|q|}\Big)^{n^2-n}(|q|)^{-2n(i+j+m)}(\bar{q})^{-(i-j)(i+j+m)}a^{-(i+j)}c^{i-j}\\
&  & \Big(\sum_{k=0}^{m+i+j}(-1)^k|q|^{k^2+k-2k(i+j+m)}{m+i+j\choose k}_{|q|^2}(bb^*)^{k+n}\Big)D^{l+j}\,.
\end{IEEEeqnarray*}
Hence, we have the following
\begin{IEEEeqnarray*}{rCl}
t_{ij}^l & = & \sum_{\substack{n=0\\ m+n=l-i\\ 0\leq m \leq l-j\\}}^{l+j}q^{n(i-j+n)}\frac{{2l\choose l+j}_{|q|^2}^{1/2}}{{2l\choose l+i}_{|q|^2}^{1/2}}{l-j\choose l-i-n}_{|q|^2}{l+j\choose n}_{|q|^2}a^mc^{l-j-m}b^nd^{l+j-n}\\
& = & a^{-(i+j)}c^{i-j}\frac{{2l\choose l+j}_{|q|^2}^{1/2}}{{2l\choose l+i}_{|q|^2}^{1/2}}\,\sum_{n=0}^{l+j}\alpha_n\,q^{n(i-j+n)}{l-j\choose l-i-n}_{|q|^2}{l+j\choose n}_{|q|^2}\\
&  & \Big(\sum_{k=0}^{l+j-n}(-1)^k|q|^{k^2+k-2k(l+j-n)}{l+j-n\choose k}_{|q|^2}(bb^*)^{k+n}\Big)D^{l+j}
\end{IEEEeqnarray*}
where, $\alpha_n=(-\bar{q})^n\Big(\frac{\bar{q}}{|q|}\Big)^{n^2-n}(|q|)^{-2n(l+j-n)}(\bar{q})^{-(i-j)(l+j-n)}$.
Now,
\begin{IEEEeqnarray*}{rCl}
&  & \sum_{n=0}^{l+j}\sum_{k=0}^{l+j-n}(-1)^k\alpha_n\,q^{n(i-j+n)}{l-j\choose l-i-n}_{|q|^2}{l+j\choose n}_{|q|^2}|q|^{k^2+k-2k(l+j-n)}{l+j-n\choose k}_{|q|^2}(bb^*)^{k+n}\\
& = & \sum_{r=0}^{l+j}\sum_{k=0}^r(-1)^k\alpha_{r-k}\,q^{(r-k)(i-j+r-k)}{l-j\choose l-i-r+k}_{|q|^2}{l+j\choose r-k}_{|q|^2}|q|^{k^2+k-2k(l+j-r+k)}{l+j-r+k\choose k}_{|q|^2}(bb^*)^r
\end{IEEEeqnarray*}
because of the fact that $\,\sum_{p=0}^s\sum_{k=0}^{s-p}f(k,p)\xi^{k+p}=\sum_{r=0}^s\Big(\sum_{k=0}^rf(k,r-k)\Big)\xi^r$. The coefficient of $(bb^*)^r$ becomes the following,
\begin{IEEEeqnarray*}{rCl}
&  & \sum_{k=0}^r(-1)^k\alpha_{r-k}\,q^{(r-k)(i-j+r-k)}|q|^{k^2+k-2k(l+j-r+k)}{l-j\choose l-i-r+k}_{|q|^2}{l+j\choose r-k}_{|q|^2}{l+j-r+k\choose k}_{|q|^2}\\
& = & \sum_{p=0}^r(-1)^{(r-p)}\alpha_p\,q^{p(i-j+p)}|q|^{(r-p)^2+(r-p)-2(r-p)(l+j-p)}{l-j\choose l-i-p}_{|q|^2}{l+j\choose p}_{|q|^2}{l+j-p\choose r-p}_{|q|^2}\\
& = & \sum_{p=0}^r(-1)^{(r-p)}\alpha_p\,q^{p(i-j+p)}|q|^{(r-p)^2+(r-p)-2(r-p)(l+j-p)}{l-j\choose l-i-p}_{|q|^2}{l+j\choose l+j-r}_{|q|^2}{r\choose r-p}_{|q|^2}\\
& = & {l+j\choose l+j-r}_{|q|^2}\,\sum_{p=0}^r\widetilde{\alpha_p}{l-j\choose l-i-p}_{|q|^2}{r\choose r-p}_{|q|^2}|q|^{2p(i-j+p)}
\end{IEEEeqnarray*}
where,
\begin{IEEEeqnarray*}{rCl}
\widetilde{\alpha_p} & = & \alpha_p(-1)^{r-p}|q|^{(r-p)^2+(r-p)-2(r-p)(l+j-p)}(\bar{q})^{-p(i-j+p)}\\
& = & (-1)^{r-p}|q|^{(r-p)^2+(r-p)-2(r-p)(l+j-p)}(\bar{q})^{-p(i-j+p)}(-\bar{q})^p|q|^{-2p(l+j-p)}\Big(\frac{\bar{q}}{|q|}\Big)^{p^2-p}(\bar{q})^{-(i-j)(l+j-p)}\\
& = & (-1)^r(\bar{q})^{(j-i)(l+j)}|q|^{r^2+r-2r(l+j)}\,\,,
\end{IEEEeqnarray*}
i,e. independent of `$p$'. Thus, the coefficient of $(bb^*)^r$ finally becomes the following,
\begin{IEEEeqnarray*}{rCl}
&  & {l+j\choose l+j-r}_{|q|^2}\,\sum_{p=0}^r\widetilde{\alpha_p}{l-j\choose l-i-p}_{|q|^2}{r\choose r-p}_{|q|^2}|q|^{2p(i-j+p)}\\
& = & {l+j\choose l+j-r}_{|q|^2}(-1)^r(\bar{q})^{(j-i)(l+j)}|q|^{r^2+r-2r(l+j)}\,\sum_{p=0}^r{l-j\choose l-i-p}_{|q|^2}{r\choose r-p}_{|q|^2}|q|^{2p(i-j+p)}\\
& = & {l+j\choose l+j-r}_{|q|^2}{l-j+r\choose l-i}_{|q|^2}(-1)^r(\bar{q})^{(j-i)(l+j)}|q|^{r^2+r-2r(l+j)}\\
& = & (\bar{q})^{(j-i)(l+j)}{l-j+r\choose i-j+r}_{|q|^2}{l+j\choose r}_{|q|^2}(-1)^r|q|^{r^2+r-2r(l+j)}\\
& = & (\bar{q})^{(j-i)(l+j)}{l-j\choose i-j}_{|q|^2}\frac{(|q|^{2(l-j+1)};|q|^2)_r}{(|q|^{2(i-j+1)};|q|^2)_r}\frac{(|q|^{-2(l+j)};|q|^2)_r}{(|q|^2;|q|^2)_r}|q|^{2r}\,\,.
\end{IEEEeqnarray*}
Therefore, the matrix coefficient $t^l_{ij}D^k$ is the following,
\begin{IEEEeqnarray*}{rCl}
t^l_{ij}D^k & = & a^{-(i+j)}c^{i-j}(\bar{q})^{(j-i)(l+j)}\frac{{2l\choose l+j}^{1/2}_{|q|^2}}{{2l\choose l+i}^{1/2}_{|q|^2}}{l-j\choose i-j}_{|q|^2}\\
&  & \Big(\sum_{r=0}^{l+j}\frac{(|q|^{-2(l+j)};|q|^2)_r}{(|q|^2;|q|^2)_r}\frac{(|q|^{2(l-j+1)};|q|^2)_r}{(|q|^{2(i-j+1)};|q|^2)_r}(|q|^2bb^*)^r\Big)D^{l+j+k}\\
& = & a^{-(i+j)}c^{i-j}(\bar{q})^{(j-i)(l+j)}\frac{{2l\choose l+j}^{1/2}_{|q|^2}}{{2l\choose l+i}^{1/2}_{|q|^2}}{l-j\choose i-j}_{|q|^2}\mathcal{P}_{l+j}^{(i-j,-i-j)}(bb^*;|q|^2)D^{l+j+k}\,.
\end{IEEEeqnarray*}
\qed

\subsection{The Fourier transform}

To each $f\in\mathcal{A}_q=\mathcal{O}(U_q(2))$, consider a matrix $\hat{f}^{(l,k)}=(\hat{f}^{(l,k)}_{m,n})\in M_{2l+1}(\mathbb{C})$ defined by $\hat{f}^{(l,k)}_{m,n}:=h(S(t^l_{m,n}D^k)f),\,l\in\frac{1}{2}\mathbb{N},\,m,n\in I_l,\,k\in\mathbb{Z}\,,$ where $h$ is the Haar state on $U_q(2)$. The mapping
\begin{align*}
\mathcal{F}:\mathcal{A}_q &\longrightarrow \widehat{\mathcal{A}_q}:=\bigoplus_{(l,k)\in\frac{1}{2}\mathbb{N}\times\mathbb{Z}}M_{2l+1,k}(\mathbb{C})\\
f &\longmapsto \hat{f}=(\hat{f}^{(l,k)})
\end{align*}
where each $M_{2l+1,k}(\mathbb{C})=M_{2l+1}(\mathbb{C})$, is called the Fourier transform on the quantum group $U_q(2)$. Let $\tau_l:M_{2l+1}(\mathbb{C})\longrightarrow\mathbb{C}$ be the $q$-trace defined by $\tau_l(M)=\sum_{i=-l}^l|q|^{-2i}m_{ii}$ for $M=(m_{ij})$. We get an inner-product on $M_{2l+1}(\mathbb{C})$ defined by $\langle M,N\rangle_l:=\tau_l(M^*N)$.
\medskip

Let $L^2(U_q(2))$ be the Hilbert space associated with $\mathcal{A}_q$ with respect to the inner-product $\langle\,x,y\rangle_L:=h(x^*y)$, and $\ell^2(\widehat{U_q(2)})$ be the Hilbert space completion of $\bigoplus_{(l,k)\in\frac{1}{2}\mathbb{N}\times\mathbb{Z}}M_{2l+1,k}(\mathbb{C})$, where each $M_{2l+1,k}(\mathbb{C})=M_{2l+1}(\mathbb{C})$, with respect to the following inner-product
\begin{center}
$\langle M,N\rangle=\sum_{l\in\frac{1}{2}\mathbb{N}}\sum_{k\in\mathbb{Z}}|2l+1|_{|q|}\langle\,M^{(l,k)},N^{(l,k)}\rangle_l\,.$
\end{center}
\bthm
The Fourier transform $\mathcal{F}:\mathcal{A}_q\longrightarrow\widehat{\mathcal{A}_q}$ is a $\mathbb{C}$-isomorphism. The inverse transform is given by
\begin{IEEEeqnarray*}{rCl}
f &=& \sum_{l\in\frac{1}{2}\mathbb{N}}\sum_{k\in\mathbb{Z}}|2l+1|_{|q|}\sum_{i,j\in I_l}|q|^{-2i}\overline{\hat{f}^{(l,k)}_{j,i}}t^l_{ij}D^k\,,
\end{IEEEeqnarray*}
and the following Plancherel formula holds
\begin{IEEEeqnarray*}{rCl}
\langle f,g\rangle_L & = & \sum_{l\in\frac{1}{2}\mathbb{N}}\sum_{k\in\mathbb{Z}}|2l+1|_{|q|}\sum_{i,j\in I_l}|q|^{-2i}\overline{\hat{f}^{(l,k)}_{j,i}}{\hat{g}^{(l,k)}_{j,i}} = \sum_{l\in\frac{1}{2}\mathbb{N}}\sum_{k\in\mathbb{Z}}|2l+1|_{|q|}\langle \hat{f}^{(l,k)},\hat{g}^{(l,k)}\rangle_l\,\,.
\end{IEEEeqnarray*}
The Fourier transform implements a unitary equivalence between the Hilbert spaces $L^2(U_q(2))$ and $\ell^2(\widehat{U_q(2)})$.
\ethm
\prf The Peter-Weyl decomposition (Thm. \ref{Peter}) shows that any $\,f\in\mathcal{A}_q$ is expressed as a finite sum of the following form
\begin{IEEEeqnarray}{rCl}\label{Fourier eq 1}
f &=& \sum_{l\in\frac{1}{2}\mathbb{N}}\sum_{k\in\mathbb{Z}}\sum_{i,j\in I_l}c^{(l)}_{i,j,k}t^{(l)}_{ij}D^k\,,
\end{IEEEeqnarray}
where the coefficients $c^{(l)}_{i,j,k}$ are given by $\langle f,t^{(l)}_{ij}D^k\rangle_L\langle t^{(l)}_{ij}D^k,t^{(l)}_{ij}D^k\rangle_L^{-1}$. Now,
\begin{IEEEeqnarray}{rCl}\label{Fourier eq 2}
\langle t^{(l)}_{ij}D^k,f\rangle_L &=& h((t^{(l)}_{ij}D^k)^*f)= h(S(t^{(l)}_{ji}D^k)f)=\hat{f}^{(l,k)}_{j,i}\,.
\end{IEEEeqnarray} 
Therefore, eqn. \ref{Fourier eq 1} and part $(ii)$ of Thm. \ref{Peter} gives us the following inverse transform,
\begin{IEEEeqnarray*}{rCl}
f &=& \sum_{l\in\frac{1}{2}\mathbb{N}}\sum_{k\in\mathbb{Z}}|2l+1|_{|q|}\sum_{i,j\in I_l}|q|^{-2i}\overline{\hat{f}^{(l,k)}_{j,i}}t^l_{ij}D^k\,.
\end{IEEEeqnarray*}
Now, for any $g\in\mathcal{A}_q$ we get that
\begin{IEEEeqnarray*}{rCl}
\langle f,g\rangle_L &=& \sum_{l\in\frac{1}{2}\mathbb{N}}\sum_{k\in\mathbb{Z}}|2l+1|_{|q|}\sum_{i,j\in I_l}|q|^{-2i}\overline{\hat{f}^{(l,k)}_{j,i}}\langle t^l_{ij}D^k,g\rangle_L\,.
\end{IEEEeqnarray*}
Since, $\langle t^l_{ij}D^k,g\rangle_L=\hat{g}^{(l,k)}_{j,i}$ by eqn. \ref{Fourier eq 2}, we get the following Plancherel formula,
\begin{IEEEeqnarray*}{rCl}
\langle f,g\rangle_L & = & \sum_{l\in\frac{1}{2}\mathbb{N}}\sum_{k\in\mathbb{Z}}|2l+1|_{|q|}\sum_{i,j\in I_l}|q|^{-2i}\overline{\hat{f}^{(l,k)}_{j,i}}\hat{g}^{(l,k)}_{j,i}\\
&=& \sum_{l\in\frac{1}{2}\mathbb{N}}\sum_{k\in\mathbb{Z}}|2l+1|_{|q|}\,\sum_i|q|^{-2i}\sum_j\overline{\hat{f}^{(l,k)}_{j,i}}\hat{g}^{(l,k)}_{j,i}\\
&=& \sum_{l\in\frac{1}{2}\mathbb{N}}\sum_{k\in\mathbb{Z}}|2l+1|_{|q|}\,\tau_l\Big((\hat{f}^{(l,k)})^*\hat{g}^{(l,k)}\Big)\\
&=& \sum_{l\in\frac{1}{2}\mathbb{N}}\sum_{k\in\mathbb{Z}}|2l+1|_{|q|}\langle \hat{f}^{(l,k)},\hat{g}^{(l,k)}\rangle_{L,l}\,\,,
\end{IEEEeqnarray*}
proving that $\mathcal{F}$ is isometry. In order to see surjectivity of $\mathcal{F}$, observe that for any $\phi=(\hat{\phi}^{(l,k)})\in\widehat{\mathcal{A}_q}$ if we choose
\begin{IEEEeqnarray*}{rCl}
f &=& \sum_{l\in\frac{1}{2}\mathbb{N}}\sum_{k\in\mathbb{Z}}\sum_{i,j\in I_l}|q|^{-2(l+j)}\frac{1-|q|^{2(2l+1)}}{1-|q|^2}\overline{\hat{\phi}^{(l,k)}_{i,j}}t^l_{j,i}D^k\in\mathcal{A}_q\,,
\end{IEEEeqnarray*}
then it follows that $\mathcal{F}(f)^{(l,k)}_{i,j}=\hat{\phi}^{(l,k)}_{i,j}\,$, and consequently $\mathcal{F}(f)=\phi\,$. This is because by eqn. \ref{Fourier eq 2}, we have
\begin{IEEEeqnarray*}{rCl}
\mathcal{F}(f)^{(l,k)}_{i,j} &=& h\Big(S(t^l_{i,j}D^k)\mathcal{F}(f)\Big) = h\Big(\Big(t^l_{j,i}D^k\Big)^*\mathcal{F}(f)\Big) = \langle t^l_{j,i}D^k,\mathcal{F}(f)\rangle_L = \widehat{\mathcal{F}(f)}^{(l,k)}_{i,j}\,.
\end{IEEEeqnarray*}
This completes the proof.\qed
\bigskip


\newsection{Tensor product decomposition}\label{Sec 5}

This section describes how the tensor product of two irreducible representations of $U_q(2)$ decomposes into irreducible 
components. In the following, we denote the character of a representation $T$ by $\chi(T)$.
\bthm\label{tensor}
One has the following,
\begin{center}
$T_{l}\otimes T_{\frac{1}{2}} \simeq T_{l+\frac{1}{2}} \oplus T_{l-\frac{1}{2}}D\,.$
\end{center}
\ethm
\prf Let 
\begin{IEEEeqnarray*}{rCl}
T_{l}\otimes T_{\frac{1}{2}} \simeq \bigoplus_{\substack{n \in \frac{1}{2}\bbn\\ m \in \bbz}} \quad C_{nm}T_nD^m.
\end{IEEEeqnarray*}
where $C_{nm}$ denotes the multiplicity of the corepresentation $T_nD^m$ in  $T_{l}\otimes T_{\frac{1}{2}}$. We have from part $(i)$ of Propn. \ref{matrixcoefficients},
\begin{IEEEeqnarray}{rCl}\label{eq3}
C_{nm}\chi(T_nD^m)&=&  \chi(T_{l}\otimes T_{\frac{1}{2}})=\chi(T_{l})\chi( T_{\frac{1}{2}})=\Big(\sum_{i=-l}^lt_{ii}^l\Big)(a+a^*D)\nonumber\\
&=&a^{2l+1}+\sum_{i=-l+1}^{l}\Big(t_{i-1,i-1}^la^*D+t_{ii}^la\Big)+t_{ll}^la^*D.
\end{IEEEeqnarray}
From Lemma \ref{class} and part $(i,iii)$ of the Propn. $\ref{clafurther}$, we have for $ -l+1\leq i \leq l$
\begin{IEEEeqnarray}{rCl}\label{eq4}
&  & t_{ll}^la^*D\,\in\cla_q[-2l-1,-2l-1,2l+1]\nonumber\\
&  & t_{i-1,i-1}^la^*D+t_{i,i}^la\,\in\cla_q[-2i+1,-2i+1,l+i] 
\end{IEEEeqnarray}
Hence, by Propn. \ref{orthogonal} and Thm \ref{faithful} we get that 
\begin{center}
$\langle \chi(T_{l}\otimes T_{\frac{1}{2}}), a^{2l+1} \rangle_L= \langle a^{2l+1} , a^{2l+1} \rangle_L>0\,.$
\end{center}
On the other hand, it follows from Thm. $\ref{Peter}$ that 
\[
 \langle \chi(T_{n}D^m), a^{2l+1} \rangle_L =\begin{cases}
                                              \langle a^{2l+1}, a^{2l+1} \rangle_L & \mbox{ if } n=l+\frac{1}{2}, m=0\,; \cr
                                              0 & \mbox{ otherwise };\cr
                                            \end{cases}
\]
since, $a^{2l+1}=t^{l+\frac{1}{2}}_{-l-\frac{1}{2},-l-\frac{1}{2}}$ by part $(i)$ of Propn. \ref{matrixcoefficients}. This proves that $C_{l+\frac{1}{2},\,0}=1$. Therefore, we have
\[
 \chi(T_{l}\otimes T_{\frac{1}{2}})-\chi(T_{l+\frac{1}{2}})=
 \sum_{\substack{n\in\frac{1}{2}\bbn\,,\,m\in\bbz\\ (n,m)\neq(l+\frac{1}{2},0)}}\,\,\chi(C_{nm}T_nD^m).
\]
Invoking eqns. $(\ref{eq3},\ref{eq4})$ and Propn. $\ref{clafurtherfurther}$, we have 
\begin{IEEEeqnarray*}{rCl}
\langle  \chi(T_{l}\otimes T_{\frac{1}{2}})-\chi(T_{l+\frac{1}{2}}), a^{2l-1}D \rangle_L &=& \langle  t_{-l,-l}^la^*D+t_{-l+1,-l+1}^la-t_{-l+\frac{1}{2},-l+\frac{1}{2}}^{l+\frac{1}{2}}\,,\,a^{2l-1}D\rangle_L\,.
\end{IEEEeqnarray*}
Using Propn. \ref{matrixcoefficients}, we have 
\begin{IEEEeqnarray*}{lCl}
&  & t_{-l,-l}^la^*D+t_{-l+1,-l+1}^la-t_{-l+\frac{1}{2},-l+\frac{1}{2}}^{l+\frac{1}{2}}\\
&=& a^{2l}a^*D+a^{2l-1}a^*Da-|q|^2{2l-1 \choose 2l-2}_{|q|^2}a^{2l-2}b^*bDa-a^{2l}a^*D+|q|^2{2l \choose 2l-1}_{|q|^2}a^{2l-1}b^*bD\\
&=& a^{2l-1}(1-|q|^2b^*b)D-|q|^2\Big(|q|^2+|q|^4+\cdots |q|^{2(2l-1)}\Big)a^{2l-1}b^*bD+|q|^2\Big(1+|q|^2+\cdots |q|^{2(2l-1)}\Big)a^{2l-1}b^*bD\\
&=& a^{2l-1}D-|q|^2\Big(1+|q|^2+|q|^4+\cdots |q|^{2(2l-1)}\Big)a^{2l-1}b^*bD+|q|^2\Big(1+|q|^2+|q|^4+\cdots |q|^{2(2l-1)}\Big)a^{2l-1}b^*bD\\
&=& a^{2l-1}D\,.
\end{IEEEeqnarray*}
Hence, $\langle\chi(T_{l}\otimes T_{\frac{1}{2}})-\chi(T_{l+\frac{1}{2}}), a^{2l-1}D \rangle_L=\langle a^{2l-1}D, a^{2l-1}D \rangle_L$. By Thm. \ref{Peter}, we have
\[
 \langle \chi(T_{n}D^m), a^{2l-1}D \rangle_L =\begin{cases}
                                              \langle a^{2l-1}D, a^{2l-1}D \rangle_L & \mbox{ if } n=l-\frac{1}{2}, m=1\,; \cr
                                              0 & \mbox{ otherwise }.\cr
                                            \end{cases}
\]
This shows that $C_{l-\frac{1}{2},1}=1$. Looking at the dimension of the corepresentations $T_l\,,\,T_{l+\frac{1}{2}}$ and $T_{l-\frac{1}{2}}D$, we get $C_{nm}=0$ for $ (n,m)\neq (l+\frac{1}{2},0)$ or $(l-\frac{1}{2},1)$.\qed

Let $P=\chi(T_{\frac{1}{2}})=a+Da^*$ and $\mathcal{B}$ be the  $C^*$-subalgebra of 
$\mathbb{A}_q$ generated by $P$ and $D$. Observe that $\mathcal{B}$ is a commutative $C^*$-algebra. Let $\alpha$ and $\beta$ be 
the roots of the quadratic polynomial $y^2-Py+D$, given by 
\[
 \alpha=\frac{P+\sqrt{P^2-4D}}{2}\quad,\quad\beta=\frac{P-\sqrt{P^2-4D}}{2}\,\,.
\]
Here, thanks to the $L^{\infty}$-functional calculus, we choose a square root of the normal operator $P^2-4D$ in the commutative von-Neumann algebra $\mathcal{B}^{\prime\prime}\subseteq C(U_q(2))^{\prime\prime}$. Therefore, $\alpha$ and $\beta$ commute with each other.
\blmma \label{charformula}
For $l \in \frac{1}{2}\bbn$, one has the following,
\begin{center}
$\chi(T_l)=\sum_{r=0}^{2l}\alpha^{2l-r}\beta^r\,.$
\end{center}
\elmma
\prf This obviously holds for $l=0$. Assume that it holds for $l\in \{0,\frac{1}{2},1,\cdots,n-\frac{1}{2}\}$. By Thm. $\ref{tensor}$, and the fact that $\chi(T_{\frac{1}{2}})=P=\alpha+\beta$ and $D=\alpha\beta$, we have 
\begin{IEEEeqnarray*}{rCl}
\chi(T_n)&=&\chi(T_{n-\frac{1}{2}})\chi(T_{\frac{1}{2}})-\chi(T_{n-1})D\\
&=& \Big(\sum_{r=0}^{2n-1}\alpha^{2n-1-r}\beta^r\Big)(\alpha+\beta)-\alpha\beta\Big(\sum_{r=0}^{2n-2}\alpha^{2n-2-r}\beta^r\Big)\\
&=& \sum_{r=0}^{2n-1}\alpha^{2n-r}\beta^r +\sum_{r=0}^{2n-1}\alpha^{2n-1-r}\beta^{r+1}- \sum_{r=0}^{2n-2}\alpha^{2n-1-r}\beta^{r+1}\\
&=& \sum_{r=0}^{2n}\alpha^{2n-r}\beta^r.
\end{IEEEeqnarray*}
\qed

\bthm \label{tensor1} The following decomposition holds,
\[
 T_{l_1}D^m\otimes T_{l_2}D^n \simeq T_{|l_1+l_2|}D^{m+n}\oplus T_{|l_1+l_2|-1}D^{m+n+1} \oplus\cdots\oplus  T_{|l_1-l_2|}D^{m+n+2\min\{l_1,l_2\}} .
\]
\ethm
\prf Without loss of generality, assume that $l_1\geq l_2$. To prove the assertion, it is enough to show that
\begin{center}
$\chi(T_{l_1}D^m\otimes T_{l_2}D^n)=\chi(T_{l_1+l_2}D^{m+n})+\chi(T_{l_1+l_2-1}D^{m+n+1})+\cdots +\chi(T_{l_1-l_2}D^{m+n+2l_2})\,.$
\end{center}
Since, $\chi(TD^r)=\chi(T)D^r$ for any corepresentation $T$ of $\cla_q$ and $r \in \bbz$, we need to show that
\begin{center}
$\chi(T_{l_1}\otimes T_{l_2})=\chi(T_{l_1+l_2})+\chi(T_{l_1+l_2-1}D)+\cdots +\chi(T_{l_1-l_2}D^{2l_2})\,.$
\end{center}
Now,
\begin{IEEEeqnarray}{lCl}\label{tensoreq1}
&  & \sum_{r=0}^{2l_2}\chi(T_{l_1+l_2-r}D^r)=\sum_{r=0}^{2l_2}\sum_{s=0}^{2l_1+2l_2-2r}\alpha^r\beta^r \alpha^{2l_1+2l_2-2r-s}\beta^{s}\nonumber\\
&=& \sum_{r=0}^{l_2}\sum_{s=0}^{2l_1+2l_2-2r}\alpha^{2l_1+2l_2-r-s}\beta^{r+s}+\sum_{r=l_2+1}^{2l_2}\sum_{s=0}^{2l_1+2l_2-2r}\alpha^{2l_1+2l_2-r-s}\beta^{r+s}\nonumber\\
&=& \sum_{r=0}^{l_2}\sum_{s=0}^{2l_1}\alpha^{2l_1+2l_2-r-s}\beta^{r+s}+\sum_{r=0}^{l_2-1}\sum_{s=2l_1+1}^{2l_1+2l_2-2r} \alpha^{2l_1+2l_2-r-s}\beta^{r+s}+\sum_{r=l_2+1}^{2l_2}\sum_{s=0}^{2l_1+2l_2-2r}\alpha^{2l_1+2l_2-r-s}\beta^{r+s}\nonumber\\
\end{IEEEeqnarray}
Let $r^{'}=2l_2-r$ and $s^{'}=2r+s-2l_2\,$. By a change of variable, we get that
\begin{IEEEeqnarray*}{lCl}
\sum_{r=0}^{l_2-1}\sum_{s=2l_1+1}^{2l_1+2l_2-2r}\alpha^{2l_1+2l_2-r-s}\beta^{r+s} &=& \sum_{r^{'}=l_2+1}^{2l_2}\,\,\sum_{s^{'}=2l_1+2l_2-2r^{'}+1}^{2l_1} \alpha^{2l_1+2l_2-r^{'}-s^{'}}\beta^{r^{'}+s^{'}}.
\end{IEEEeqnarray*}
Putting this in eqn. \ref{tensoreq1}, we finally have the following
\begin{IEEEeqnarray*}{lCl}
&  & \sum_{r=0}^{2l_2}\chi(T_{l_1+l_2-r}D^r)\\
&=& \sum_{r=0}^{l_2}\sum_{s=0}^{2l_1}\alpha^{2l_1+2l_2-r-s}\beta^{r+s}+ \sum_{r=l_2+1}^{2l_2}\sum_{s=0}^{2l_1} \alpha^{2l_1+2l_2-r-s}\beta^{r+s}\\
&=& \sum_{r=0}^{2l_2}\sum_{s=0}^{2l_1}\alpha^{2l_1+2l_2-r-s}\beta^{r+s}\\
&=& \sum_{r=0}^{2l_2} \alpha^{2l_2-r}\beta^r\sum_{s=0}^{2l_1}\alpha^{2l_1-s}\beta^{s}\\
&=& \chi(T_{l_1}\otimes T_{l_2})\,,
\end{IEEEeqnarray*}
and this completes the proof.\qed


\newsection{The case of \texorpdfstring{$|q|=1$}{}}\label{Sec 6}

This section deals with the case of $|q|=1$ and $q$ not a root of unity. Let $\clk$ be the Hilbert space $\ell^2(\bbz)\otimes\ell^2(\bbz)\otimes\ell^2(\bbz)$ and $U:e_n\longmapsto e_{n+1}$ be the bilateral shift acting on $\ell^2(\bbz)$. Define 
the representation $\pi$ of the $C^*$-algebra $C(U_q(2))$ on $\clk$ as follows~:
\begin{IEEEeqnarray}{rCl}\label{representation1}
\pi(a) &=& U\otimes I\otimes\frac{U^*+U}{2}\,,\nonumber\\
\pi(b) &=& q^N\otimes U\otimes\frac{U^*-U}{2i}\,,\nonumber\\
\pi(D) & =& U^2\otimes I\otimes U.
\end{IEEEeqnarray}
One can view  $C(U_q(2))$ as a $C^*$-subalgebra of $\mathcal{B}\big(\ell^2(\bbz)\otimes\ell^2(\bbz)\big)\otimes C(\bbbt)$ by identifying $U$ with the function $\bbt:t\mapsto t$. For $t \in \bbbt$, let $\psi_t:\mathcal{B}\big(\ell^2(\bbz)\otimes\ell^2(\bbz)\big)\otimes C(\bbbt)\longrightarrow\mathcal{B}\big(\ell^2(\bbz)\otimes\ell^2(\bbz)\big)$ be the homomorphism given by the restriction of $\,I\otimes I\otimes ev_t$ to $C(U_q(2))$, where $ev_t$ is the evaluation at $t$. 
\bppsn\label{rep1}
The representation $\pi$ of $C(U_q(2))$ defined above is faithful.
\eppsn
\prf It is enough to show that all irreducible representations of $C(U_q(2))$ factor through the representation $\pi$. Let $\rho$ be an irreducible representation of $C(U_q(2))$. Using the description given in (Thm. $3.4$, \cite{ZhaZha-2005aa}), it is not difficult to establish that if $\rho(b)=0$ then $\rho$ factors through $\psi_1 \circ \pi$, if $\rho(a)=0$ then $\rho$ factors through $\psi_{\sqrt{-1}}\circ\pi$, and if $\rho(a),\rho(b) \neq 0$ then $\rho$ factors through $\psi_t \circ \pi$ for some $t \neq 1, \sqrt{-1}$.\qed
\bthm[\cite{ZhaZha-2005aa}]\label{haar1} 
The Haar state $h:C(U_q(2)) \longrightarrow \bbc$ is given as follows~:
\begin{IEEEeqnarray*}{rCl}
 h( \langle n,m,k,l \rangle)= \begin{cases}
                               \frac{1}{m+1} & \mbox{ if } m=k, \mbox{ and } n=l=0, \cr
                               0 & \mbox{ otherwise .}
                              \end{cases}
\end{IEEEeqnarray*}
\ethm
It is known that in this case the Haar state is trace \cite{ZhaZha-2005aa}. Consider the homomorphism $\,\phi:C(U_q(2))\longrightarrow C(\bbbt)$ given by $\phi(a)=\textbf{z},\,\phi(b)=0$, and $\phi(D)=1$. This gives a $\bbbt$-coaction $\Phi:C(U_q(2))\longrightarrow C(U_q(2))\otimes C(\bbbt)$ on $C(U_q(2))$ defined by $\Phi(x)=(\mbox{id}\otimes \phi)\circ\Delta\,$. The quotient space $\bbbt\diagdown U_q(2)\diagup\bbbt$ is defined as follows~:
\[
 C(\bbbt\diagdown U_q(2)\diagup\bbbt)=\{x \in C(U_q(2)): (\phi \otimes \mbox{id}  \otimes \phi)(\Delta \otimes \mbox{id} )  \Delta(x)=1 \otimes x\otimes 1 
 \}.
\]
In such a case, the conditional expectation $E: C(U_q(2))\longrightarrow C(\bbbt\diagdown U_q(2)\diagup\bbbt)$ is defined to be the map $((h_{\bbbt}\circ\phi)\otimes\mbox{id}\otimes(h_{\bbbt}\circ\phi))(\Delta \otimes \mbox{id})\circ\Delta\,$.
\blmma The $C^*$-algebra $C(\bbbt\diagdown U_q(2)\diagup\bbbt)$ is the $C^*$-subalgebra  of $C(U_q(2))$ generated by $D$.
\elmma
\prf Let $\mathcal{O}(U_q(2))$ be the $\star$-subalgebra of $C(U_q(2))$ generated by $a,b$ and $D$. Then, we have
\begin{IEEEeqnarray}{rCl}  \label{eq61}
C(\bbbt\diagdown U_q(2)\diagup\bbbt)= E(C(U_q(2))=E(\overline{\mathcal{O}(U_q(2))})=\overline{E(\mathcal{O}(U_q(2)))}\,.
\end{IEEEeqnarray}
One has the following,
\begin{IEEEeqnarray*}{rCl}
 (h_{\bbbt} \circ \phi)(t_{i,j}^lD^m)=\begin{cases}
                                    1 & \mbox{ if } l=0, \cr
                                    0 & \mbox{ if } l\neq 0.
                                   \end{cases}
\end{IEEEeqnarray*}
Hence,
\begin{IEEEeqnarray*}{rCl}
E(t_{i,j}^lD^m)=\sum_{r,s\in\{-l,\ldots,l\}}(h_{\bbbt} \circ \phi)(t_{i,r}^lD^m)  t_{r,s}^lD^m  (h_{\bbbt} \circ \phi)(t_{s,j}^lD^m)
  =\begin{cases}
    D^m & \mbox{ if } l=0, \cr
    0 & \mbox{ if } l\neq 0.
   \end{cases}
\end{IEEEeqnarray*}
Using this in eqn. \ref{eq61} we get our claim.\qed
\bthm\label{faithful1}
The Haar state $h$ on the quantum group $U_q(2)$ is faithful.
\ethm
\prf Let $\mathcal{C}:=C^*(D)$ be the $C^*$-subalgebra of $C(U_q(2))$ generated by $D$. Invoking the argument used in Thm. $\ref{faithful}$, it is enough to show that $\left.h\right|_\mathcal{C}$ is faithful. Observe that  the map $D\mapsto\textbf{z}$ is a $C^*$-algebra isomorphism between $\mathcal{C}$ and $C(\bbbt)$. For $m\in\bbz$, we have by Thm. $\ref{haar1}$
\[
 \left.h\right|_\mathcal{C}(D^m)=\begin{cases}
              1 & \mbox{ if } m=0, \cr
              0 & \mbox{ otherwise }.
             \end{cases}
\]
Hence, $\left.h\right|_\mathcal{C}$ is same as the Haar state on $C(\bbbt)$, which is faithful.\qed

With Propn. \ref{rep1} and Thm. \ref{faithful1} in hand, one can check that all the results obtained in Sec \S \ref{Sec 3} holds in this case also. Few results obtained in Sec \S \ref{Sec 4} change a bit and here we indicate those without proof, as these can be proved similarly.
\blmma \label{induction2}
For $r,s\in\bbn$, let $d_{r,s}=\|a^rb^s\|_L^2$ and $c_{r,s}=\|a^rb^s\|_R^2\,$. Then, one has the followings.
\begin{enumerate}
\item[(i)] $d_{r,s}=d_{r-1,s}-d_{r-1,s+1}\,\,\forall\,r\geq 1\,;$
\item[(ii)] $d_{r,s}=\Big((r+s+1){ r+s \choose s}\Big)^{-1}\,;$
\item[(iii)] $d_{r,s}=c_{r,s}\,$.
\end{enumerate}
\elmma
\prf The last part follows from the traciality of the Haar state.\qed

The Peter-Weyl decomposition in this case is obtained in \cite{Zha-2010aa}, but here we get the following betterment including the norm factor of the matrix coefficients.
\bthm
\begin{enumerate}[(i)]
\item For $l\in\frac{1}{2}\bbn$, let $t_{ij}^l,\,i,j\in I_l$ be the matrix coefficients of the corepresentation $T_l^R$ with respect to the basis $\{f_j^l: j \in I_l\}$. Then, we have
\begin{IEEEeqnarray*}{rCl}
t_{ij}^l= \sum_{\substack{m+n=l-i\\ 0\leq m \leq l-j\\0\leq n \leq l+j}}q^{n(l-j-m)}
\frac{{ 2l \choose l+j}^{1/2}}{{ 2l \choose l+i}^{1/2}}{l-j \choose m}{l+j \choose n}a^mc^{l-j-m}b^nd^{l+j-n}\,;
\end{IEEEeqnarray*}
\item $\langle t_{ij}^{l}D^m, t_{i^{'}j^{'}}^{l^{'}}D^{m^{'}}\rangle=(2l+1)^{-1}\delta_{ll^{'}}\delta_{ii^{'}}\delta_{jj^{'}}\delta_{mm^{'}}\,;$
\item The set $\big\{\sqrt{2l+1}\,t_{ij}^{l}D^m:l \in \frac{1}{2}\bbn, m \in \bbz\big\}$ is an orthonormal basis of $L^2(h)$.
\end{enumerate}
\ethm
\bthm
The Fourier transform $\mathcal{F}:\mathcal{A}_q\longrightarrow\widehat{\mathcal{A}_q}$ is a $\mathbb{C}$-isomorphism. The inverse transform is given by
\begin{IEEEeqnarray*}{rCl}
f &=& \sum_{l\in\frac{1}{2}\mathbb{N}}\sum_{k\in\mathbb{Z}}(2l+1)\sum_{i,j\in I_l}\overline{\hat{f}^{(l,k)}_{j,i}}t^l_{i,j}D^k\,,
\end{IEEEeqnarray*}
and the following Plancherel formula
\begin{IEEEeqnarray*}{rCl}
\langle f,g\rangle & = & \sum_{l\in\frac{1}{2}\mathbb{N}}\sum_{k\in\mathbb{Z}}(2l+1)\sum_{i,j\in I_l}\overline{\hat{f}^{(l,k)}_{j,i}}{\hat{g}^{(l,k)}_{j,i}} = \sum_{l\in\frac{1}{2}\mathbb{N}}\sum_{k\in\mathbb{Z}}(2l+1)\langle \hat{f}^{(l,k)},\hat{g}^{(l,k)}\rangle_l
\end{IEEEeqnarray*}
holds. Here, the notations are same as defined in subsection $4.4$, but the $q$-trace involved in $\langle\,.\,,\,.\rangle_l$ becomes the usual matrix trace.
\ethm
\bigskip


\section{Classification of \texorpdfstring{$U_q(2)$}{}}\label{Sec 7}

In this section, we classify the compact quantum group $U_q(2)$ for $q\in\mathbb{C}^*$ and not roots of unity. This subsumes earlier investigation in \cite{Zha2-2010aa}. For the case of $SU_q(2),\,q\in\mathbb{R}^*$, see \cite{Wang-1999aa}. Since, we will be dealing with different values of $q$ at the same time, the generators of $U_q(2)$ are denoted by $a_q,b_q$ and $D_q$ in order to avoid any confusion. Let $\Omega$ denotes the set of all roots of unity.
 
\bthm \label{classification}
Let $q$ and $q^{'}$ be two non-zero complex numbers which are not roots of unity. 
Then, $U_q(2)$ and $U_{q^{'}}(2)$ are isomorphic if and only $q^{'}\in\{q,\overline{q},\frac{1}{q},\frac{1}{\overline{q}}\}$.
\ethm
\prf \textbf{Step 1:} For $q\in\mathbb{C}^*\setminus\Omega,\,U_q(2)$ and $U_{\frac{1}{q}}(2)$ are isomorphic.\\
Let $q^{'}=1/q$. Define the map $\Phi:C(U_q(2))\longrightarrow C(U_{q^{'}}(2))$ by 
$\Phi(a_q)=a_{q^{'}}^*D_{q^{'}},\,\Phi(b_q)=-q^{'}b_{q^{'}}^*D_{q^{'}},\,\Phi(D_q)=D_{q^{'}}$. Using the relations in \ref{relations}, one can show  that 
$\Phi$ is a  $C^*$-algebra homomorphism. The inverse map $\Psi:C(U_{q^{'}}(2))\longrightarrow C(U_q(2))$ is given by $a_{q^{'}}\mapsto a_qD_q^*,\, b_{q^{'}}\mapsto -qb_q^*D_q,\, D_{q^{'}}\mapsto D_q$. Moreover, to show that $\Delta_{q^{'}}\circ\Phi=(\Phi\otimes\Phi)\circ\Delta_{q}$ it is enough to check this condition on the generators $a_q,b_q$ and $D_q$  of $C(U_q(2))$, since $\Phi$ is a homomorphism. This follows easily from the eqn. \ref{comul}.
\medskip

\textbf{Step 2:} For $q\in\mathbb{C}^*\setminus\Omega,\,U_q(2)$ and $U_{\frac{1}{\overline{q}}}(2)$ are isomorphic.\\
Let $q^{'}=1/\overline{q}$. Define the map $\Phi:C(U_q(2))\longrightarrow C(U_{q^{'}}(2))$ by $\Phi(a_q)=a_{q^{'}}^*,\,\Phi(b_q)=\overline{q^{'}}b_{q^{'}}^*,\,\Phi(D_q)=D_{q^{'}}^*$. Using the relations in \ref{relations}, one can show that $\Phi$ is a $C^*$-algebra homomorphism. The inverse map $\Psi:C(U_{q^{'}}(2))\longrightarrow C(U_q(2))$ is given by $a_{q^{'}}\mapsto a_q^*,\, b_{q^{'}}\mapsto \overline{q}b_q^*,\, D_{q^{'}}\mapsto D_q^*$. Moreover, to show that $\Delta_{q^{'}}\circ\Phi=(\Phi\otimes\Phi)\circ\Delta_{q}$ it is enough to check this condition on the generators $a_q,b_q$ and $D_q$  of $C(U_q(2))$, since $\Phi$ is a homomorphism. This follows easily from the eqn. \ref{comul}.
\medskip

\textbf{Step 3:} For $q\in\mathbb{C}^*\setminus\Omega,\,U_q(2)$ and $U_{\overline{q}}(2)$ are isomorphic.\\
This follows from Step $1$ and $2$.
\medskip

\textbf{Step 4:} For $q,q^{'}\in\bbc^*\setminus\Omega$, let $U_q(2)$ and $U_{q^{'}}(2)$ be isomorphic. Then, $q^{'}\in\{q,\overline{q},\frac{1}{q},\frac{1}{\overline{q}}\}$.\\
Let $\varPsi:C(U_q(2))\longrightarrow C(U_{q^{'}}(2))$ be a quantum group isomorphism.
Then, $\varPsi$ induces a bijective correspondence between one dimensional representations of $U_q(2)$ and $U_{q^{'}}(2)$. From Thm. \ref{Peter}, for one dimensional representations one has $l=0$, and thus it is $D^k$ for $k\in\bbz$. Hence, $\varPsi(D_q)=D_{q^{'}}^k$ with $k\in\bbz$. Since $\varPsi$ is bijective, one has $k=\pm 1$ i,e., either $\varPsi(D_q)=D_{q^{'}}$ or $\varPsi(D_q)=D_{q^{'}}^*$. First assume that $\varPsi(D_q)=D_{q^{'}}$. Let 
\begin{IEEEeqnarray*}{rCl}
u=\left[ {\begin{matrix}
   a_q  & b_q\\
   -qb_q^* D_q& D_qa_q^* \\
  \end{matrix}} \right]\quad,\quad v_m= \left[ {\begin{matrix}
   a_{q^{'}}D_{q^{'}}^m  & b_{q^{'}}D_{q^{'}}^m \\
   -q^{'}b_{q^{'}}^*D_{q^{'}}^{m+1}  & a_{q^{'}}^*D_{q^{'}}^{m+1}  \\
  \end{matrix}} \right].
\end{IEEEeqnarray*}
It is known that two equivalent unitary representations $\pi_1,\pi_2$ of a compact quantum group are unitarily equivalent. Since, $\varPsi(u)$ is unitarily equivalent to $v_m$ for some $m \in \bbz$, there exists a $2 \times 2$ unitary matrix $C=\big(c_{ij}\big)$ such that
\begin{center}
$v_m=C\varPsi(u)C^{-1}=\varPsi(CuC^{-1})\,.$
\end{center}
Then, we have
\begin{IEEEeqnarray}{lCl}\label{used here}
(CuC^{-1})_{11}(D^*_q)^m &=& \Big((CuC^{-1})_{22}(D^*_q)^{m+1}\Big)^*\nonumber\\
(CuC^{-1})_{12}(D_q^*)^{m} &=& \Big(-\frac{1}{q^{'}}(CuC^{-1})_{21}(D_q^*)^{m+1}\Big)^*
\end{IEEEeqnarray}
This follows from the observation that applying $\varPsi$ on both sides they agree, and $\varPsi$ is an isomorphism.
Using the fact that
\begin{center}
$C^{-1}=C^*= \left[ {\begin{matrix}
   \overline{c_{11}}  & \overline{c_{21}}\\
  \overline{c_{12}} & \overline{c_{22}} \\
  \end{matrix} } \right]\,,$
  \end{center}
from eqn. $\ref{used here}$ we get the following equations,
\begin{IEEEeqnarray}{lCl}\label{claeq1}
&  & c_{11}\overline{c_{11}}a_q(D_q^*)^m-qc_{12}\overline{c_{11}}b_q^*(D_q^*)^{m-1}+c_{11}\overline{c_{12}}b_q(D_q^*)^m+c_{12}\overline{c_{12}}a_q^*(D_q^*)^{m-1}\nonumber\\
&=& c_{21}\overline{c_{21}}a_q^*D_q^{m+1}-\bar{q}\Big(\frac{|q|}{q}\Big)^{2m}c_{21}\overline{c_{22}}b_qD_q^{m}+\Big(\frac{|q|}{\overline{q}}\Big)^{2(m+1)}\overline{c_{21}}c_{22}b_q^*D_q^{m+1}+\overline{c_{22}}c_{22}a_qD_q^{m}
\end{IEEEeqnarray}
and
\begin{IEEEeqnarray}{lCl}\label{claeq2}
&  & c_{11}\overline{c_{21}}a_q(D_q^*)^m-qc_{12}\overline{c_{21}}b_q^*(D_q^*)^{m-1}+c_{11}\overline{c_{22}} b_q(D_q^*)^m +c_{12}\overline{c_{22}}a_q^*(D_q^*)^{m-1}\nonumber \\
&=& -\frac{1}{\overline{q^{'}}}\overline{c_{21}}c_{11}a_q^*D_q^{m+1}+\frac{\bar{q}}{\overline{q^{'}}}\Big(\frac{|q|}{q}\Big)^{2m}\overline{c_{22}}c_{11}b_qD_q^{m}-\frac{1}{\overline{q^{'}}}\Big(\frac{|q|}{\overline{q}}\Big)^{2(m+1)}\overline{c_{21}}c_{12}b_q^*D_q^{m+1}-\frac{1}{\overline{q^{'}}}\overline{c_{22}}c_{12}a_qD_q^{m}\,.\nonumber\\
\end{IEEEeqnarray}
Since $C=(c_{ij})$ is a non-zero matrix, i,e. not all its entries are qual to zero simultaneously, it follows from these equations that $m=0$ due to Thm. \ref{basis}. Now, equating coefficients from both sides of eqns. $(\ref{claeq1},\ref{claeq2})$, we get the following equation,
\begin{IEEEeqnarray}{lCl}\label{claeq3}
c_{11}\overline{c_{11}}=c_{22}\overline{c_{22}}\,\,,\,\,c_{12}\overline{c_{12}}=c_{21}\overline{c_{21}}\,\,,\,\,c_{11}\overline{c_{22}}=\frac{\bar{q}}{\overline{q^{'}}}\overline{c_{22}}c_{11}\,\,,\,\,\bar{q}c_{12}\overline{c_{21}}=  \frac{1}{\overline{q^{'}}}\overline{c_{21}}c_{12}\,.
\end{IEEEeqnarray}
If $c_{11}$ and $c_{12}$ both are non-zero then eqn. $\ref{claeq3}$ implies that none of the entries of the matrix $C$ are zero. This forces $q=q^{'}$ and $q=1/q^{'}$ simultaneously and hence, $q=q^{'}=\pm 1\in\Omega$. Therefore, either $c_{11}=0$ or $c_{12}=0$, as both $c_{11}$ and $c_{12}$ can not be simultaneously zero since $C$ is a non-zero matrix (again from eqn. $\ref{claeq3}$). If $c_{11}=0$ then, $c_{22}=0,\,c_{12}\neq 0$ and $c_{21}\neq 0$, and we have $q=1/q^{'}$ from eqn. $\ref{claeq3}$. If $c_{12}=0$ then, $c_{21}=0,\,c_{11}\neq 0$ and $c_{22}\neq 0$, and we have $q=q^{'}$ from eqn. $\ref{claeq3}$.\\
Now, if $\varPsi(D_q)=D_{q^{'}}^*$ then, we replace $D_q^*$ by $D_q$ in eqn. \ref{used here} and proceeding similarly as above we get the following equation,
\begin{IEEEeqnarray}{lCl}\label{claeq5}
c_{11}\overline{c_{11}}=c_{22}\overline{c_{22}}\,\,,\,\,c_{12}\overline{c_{12}}=c_{21}\overline{c_{21}}\,\,,\,\,c_{11}\overline{c_{22}}=\frac{q}{\overline{q^{'}}}\overline{c_{22}}c_{11}\,\,,\,\,qc_{12}\overline{c_{21}}=  \frac{1}{\overline{q^{'}}}\overline{c_{21}}c_{12}\,.
\end{IEEEeqnarray}
From this equation it will follow that either $q=\overline{q^{'}}$ or $q=\frac{1}{\overline{q^{'}}}$ by the same justification as above.

Combining all these steps, we get the assertion.\qed
\brmrk
This theorem justifies that for $|q|\neq 1$, it is enough to do the computations for $|q|<1$ only.
\ermrk
\bigskip

\section*{Acknowledgements}
Satyajit Guin acknowledges useful discussion with Sutanu Roy on the article \cite{KasMeyRoyWor-2016aa} and hospitality provided by NISER Bhubaneswar during his visit. He also acknowledges the support of DST INSPIRE Faculty award grant DST/INSPIRE/04/2015/000901. Bipul Saurabh acknowledges the support of SERB grant SRG/2020/000252.

\bigskip

\bigskip

\noindent{\sc Satyajit Guin} (\texttt{sguin@iitk.ac.in})\\
         {\footnotesize Department of Mathematics and Statistics,\\
         Indian Institute of Technology, Kanpur,\\
         Uttar Pradesh 208016, India}
\bigskip

\noindent{\sc Bipul Saurabh} (\texttt{bipul.saurabh@iitgn.ac.in})\\
         {\footnotesize Indian Institute of Technology, Gandhinagar,\\  Palaj, Gandhinagar 382355, India}

\end{document}